\documentclass{article}

\usepackage{latexsym}
\usepackage{amsfonts}
\usepackage[mathscr]{euscript}
\usepackage{amssymb}
\usepackage{amsthm}
\usepackage{amsmath}
\usepackage[pdftex]{graphicx,color}

\def\bb{{\cal B}}

\def\ll{{\cal L}}
\def\cc{{\cal C}}
\def\mm{{\cal M}}
\def\nn{{\cal N}}

\def\bp{{\mathbb P}}
\def\be{{\mathbb E}}

\def\og1{{\overline{G_1}}}
\def\og2{{\overline{G_2}}}

\newtheorem{theorem}{Theorem}
\newtheorem{corollary}[theorem]{Corollary}

\newtheorem{proposition}[theorem]{Proposition}
\newtheorem{lemma}[theorem]{Lemma}
\newtheorem{observation}[theorem]{Observation}

\begin{document}
\begin{center}

\vglue -0.9 cm

{\LARGE Simple Euclidean arrangements with one $(\ge 5)$--gon}

\vglue 0.1 cm

\end{center}

\vglue 0.4 cm 

\begin{center}
{\large J. Lea\~nos\footnote{Campus Jalpa, Universidad
                   Aut\'onoma de Zacatecas, Zac.,    98600 M\'exico \\ E--mail: {jesus.leanos@gmail.com}} \hglue 1.4 cm C. Ndjatchi Mbe Koua\footnote{Coordinaci\'on de Ingenier\'ia  Mecatr\'onica y Energ\'ia, Universidad Polit\'ecnica de Zacatecas (UPZ), Fresnillo, Zac.,  99000 M\'exico }
\hglue 1.4 cm L. M. Rivera--Mart\'inez$^{1,}$}\footnote{Rivera--Mart\'inez was supported by PROMEP (SEP), grant UAZ-PTC-103}
\end{center}

\vglue 0.4 cm 

\begin{center} {\large August 24, 2010} \end{center}

\vglue 0.5 cm

{\leftskip=30pt \rightskip=30pt \small \noindent{\bf Abstract.} Let $\ll$ be a simple Euclidean arrangement of $n$ pseudolines.  It is shown that if $\ll$ has exactly one $(\ge 5)$--gon $P$,  and $k$ is the number of edges of $P$ that are adjacent to an unbounded cell of the subarrangement of $\ll$ induced by the pseudolines in 
$P$, then  $\ll$ has exactly $n-k$ triangles and $k+n(n-5)/2$ quadrilaterals. We also prove that if each pseudoline of $\ll$ is adjacent to $P$ then $\ll$ is stretchable.}

\section{Introduction}

Recently, Lea\~nos et al.~\cite{leanos} proved that if a simple Euclidean arrangement of $n$ pseudolines has no $(\ge 5)$--gons, then it is stretchable and has exactly $n-2$ triangles and $(n-2)(n-3)/2$ quadrilaterals. Our goal in this paper is to study the same problems for the case when  all the bounded cells of the arrangement under consideration are,  except one,  either  triangles or  quadrilaterals.

  We recall that a simple noncontractible closed curve in the projective plane $\bp$ is a {\it pseudoline}, 
and an {\it arrangement of pseudolines} is a collection $\bb = \{p_0 , p_1, \ldots ,p_n \}$ of pseudolines that intersect (necessarily cross) pairwise exactly once. Since $\bp \setminus p_0$ is homeomorphic to the Euclidean plane  $\be$, we may regard $\{p_1, \ldots ,p_n\}$ as an {\it arrangement of pseudolines} in $\be$ (and regard $p_1,\ldots ,p_n$ as  {\it pseudolines} in $\be$). An arrangement is {\it simple} if no point belongs to more than two pseudolines.  The cell complex of an Euclidean arrangement in $\bp$ has both bounded and unbounded cells. As in~\cite{felsner-kriegel} and~\cite{leanos}, we are only interested in bounded cells (whose interiors are the  {\it polygons} or  {\it faces}). Thus it is clear what  is meant by a {\it triangle}, a {\it quadrilateral}, or, in general, an $n$--{\it gon} of the arrangement. In this work we  are interested in arrangements in which every bounded cell is, except one,  either a triangle or a quadrilateral;  for this reason we say that such an arrangement has one  $(\ge 5)$-gon.

One of the most interesting and widely studied problems concerning arrangements of lines and 
pseudolines is the determination of upper and lower bounds for the number $p_k$ of $k$--gons. An extensive amount of research in such problems fo\-llo\-wed Gr\"{u}nbaum's seminal work

\begin {theorem}\label{igual-ss} 
Let $\ll$ be a simple Euclidean arrangement of $n$ pseudolines with exactly one $(\ge 5)$--gon $P$, and  let  $k$ be the number of edges of $P$ that are adjacent to an unbounded cell of the subarrangement of $\ll$ induced by the pseudolines in  $P$. Then  $\ll$ has exactly $n-k$ triangles.
\end {theorem}

As an immediate consequence of Theorem~\ref{igual-ss} and the fact that the number of bounded cells in a simple Euclidean arrangement of $n$ pseudolines is $1+n(n-3)/2$, we have the following corollary.

\begin {corollary}\label{coro-s} 
Let $\ll$ be a simple Euclidean arrangement of $n$ pseudolines with exactly one $(\ge 5)$--gon $P$, and  let  $k$ be the number of edges of $P$ that are adjacent to an unbounded cell of the subarrangement of $\ll$ induced by the pseudolines in  $P$. Then  $\ll$ has exactly $n-k$ triangles, $k + n(n-5)/2$  quadrilaterals, and one $(\ge 5)$--gon.
\end {corollary}
\noindent{\it Proof.} Let $p_3$ and $p_4$ be the number of triangles and quadrilaterals of $\ll$, respectively. Since the number of bounded cells in a simple Euclidean arrangement of $n$ pseudolines is $1+n(n-3)/2$ and all the bounded cells of $\ll$ are, except one, either  triangles or quadrilaterals, then  $p_3+p_4 = n(n-3)/2$. By Theorem~\ref{igual-ss}, $p_3=n-k$. This implies $p_4=k + n(n-5)/2$, as desired. 
$\square$
\\

\begin{figure}[ht]
  \begin{center}
\includegraphics[width= 13 cm, height=3.9cm]{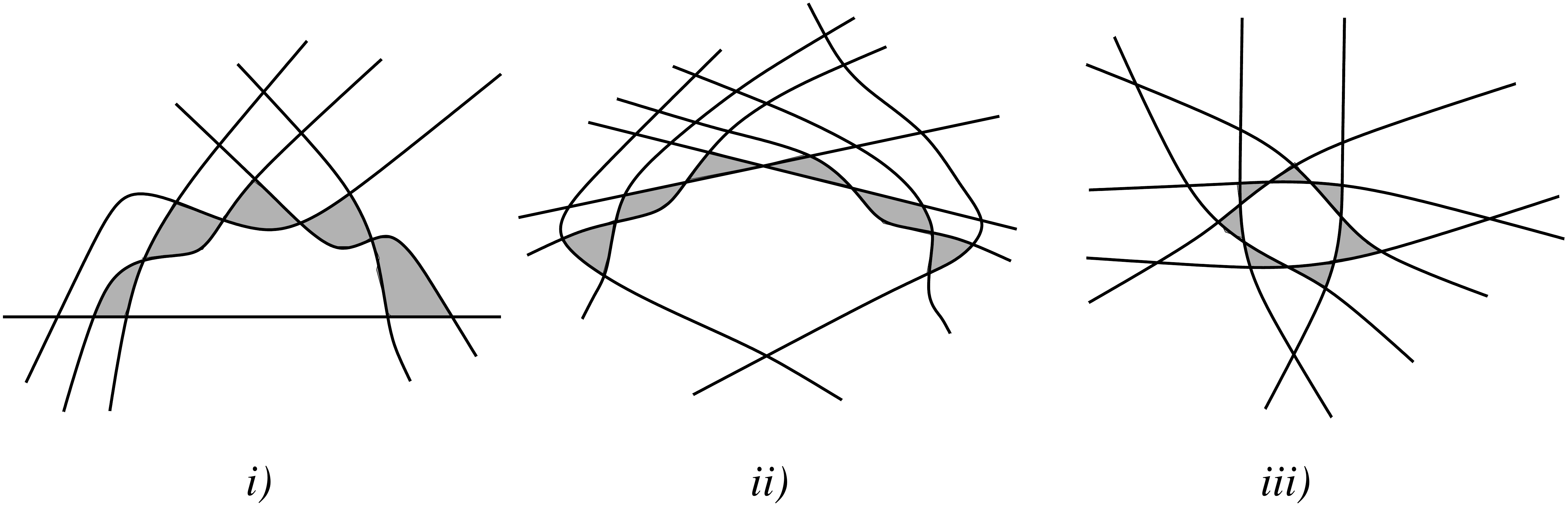}
  \end{center}
  \caption{\small Three arrangements of $\Im$.}
  \label{fig:im}
  \end{figure}

In this work, we denote by  $\Im$ the set of  simple Euclidean arrangements of  pseudolines with  one $(\ge 5)$--gon in which every pseudoline is adjacent to the (unique) $(\ge 5)$--gon. Every  arrangement shown in Figure~\ref{fig:im} is an element of $\Im$.
  
Our first aim in this work is to prove  Theorem~\ref{igual-ss} for the case where the arran\-ge\-ment in question belongs to  $\Im$ (namely, Theorem~\ref{igual-t}). This is the content of Section~\ref{somelemas}. 
In Section~\ref{general} we show Theorem~\ref{igual-ss} by induction on the number of pseudolines of $\ll$ that are not in $P$ (Theorem~\ref{igual-t} is the  base case of the induction in the proof of Theorem~\ref{igual-ss}). 

Although the main aim of all the  lemmas and propositions in Sections~\ref{somelemas} and~\ref{general} is to prove Theorem~\ref{igual-ss}, some of them are  interesting in themselves. 
 This is the case of  Proposition~\ref{coro-triangulo-adyacente}, which says that every pseudoline in a simple Euclidean arrangement is adjacent to at least one triangle.

Propositions~\ref{almenos-t} and~\ref{alomas-t} (in Section~\ref{somelemas})  
have applications not only in the proof  of Theorem~\ref{igual-ss}  but also in the proof of  our main result on stretchability (Theorem~\ref{rectif}). Together,  these propositions reveal a structural property of the arrangements of $\Im$:  {\it if $Q$ is the $(\ge 5)$--gon of $\mm\in \Im$  and $T$ is  a triangle of $\mm$, then $Q$ and $T$ have a common edge};  see the examples in Figure~\ref{fig:im}. In the proof of  Theorem~\ref{rectif} (Section~\ref{secrec}), we use this structural property to guarantee  that every arrangement in $\Im$  has $3$ pseudolines with certain properties. The properties of such $3$ pseudolines play a crucial role in the proof of Theorem~\ref{rectif}.  

    An {\it arrangement of lines} in $\be$ is a collection of straight lines, no two of them parallel. Thus, 
every arrangement of lines is an arrangement of pseudolines. On the other hand, not every 
    arrangement of pseudolines is {\it stretchable}, that is, equivalent to an arrangement of lines, where two arrangements are {\it equivalent} if they {ge\-ne\-ra\-te} isomorphic cell decompositions of $\be$. 
    Every arrangement of eight  pseudolines is stretchable~\cite{goodman-pollack}, but there is a simple non--stretchable arrangement  of nine pseudolines~\cite{ringel} (unique up to isomorphism; see~\cite{goodman-pollack}). Stretchability questions are typically 
difficult: deciding stretchability is NP--hard~\cite{shor}, even for simple arrangements~\cite{bose}. 
    The concept of stretchability is particularly relevant because of the close {co\-nnec\-tion} between 
arrangements of pseudolines and rank 3 oriented matroids: on this ground, the problem of 
    stretchability of arrangements is equivalent to the problem of realizability for oriented matroids
     (see ~\cite{bjorner, richter}).  
     
 The following theorem is the main result of this work with respect to stretchability and  establishes that every arrangement of $\Im$ is stretchable. 
 
\begin {theorem}\label{rectif} 
If  $\ll$ is an arrangement of  $\Im$, then $\ll$ is stretchable.
\end {theorem}

The proof of Theorem~\ref{rectif} will be given in Section~\ref{secrec}.  Finally, in Section~\ref{cuantos}, we show that  there are exponentially many non--isomorphic arrangements of $\Im$.


\subsection{Notation}

An {\it $n$--gonal region} $Q$ of an Euclidean arrangement of pseudolines $\ll$ is any $n$--gon of a  subarrangement of $\ll$ (thus, every $n$--gon of $\ll$ is an $n$--gonal  region of $\ll$ but not conversely). 

\begin{figure}[ht]
  \begin{center}
  \includegraphics[width= 12 cm, height=4.5cm]{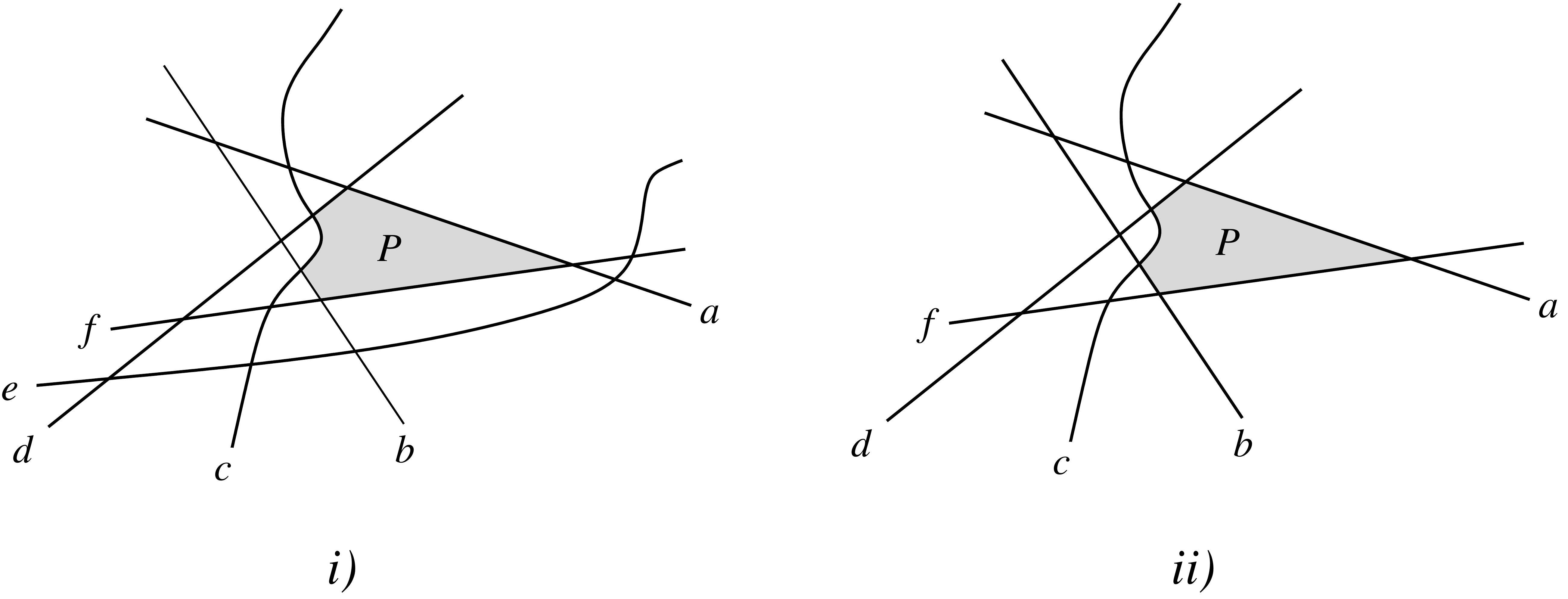}
  \end{center}
  \caption{\small Let $\ll:=\{a,b,c,d, e, f\}$ be the arrangement in $i)$ and let  $\ll'$ be the arrangement in $ii)$. Note that $P$ is the only $(\ge 5)$--gon of $\ll$, and that $\ll'$ is the subarrangement of  $\ll$  induced by the pseudolines in $P$. Since $e_{a,P}$ and $e_{f,P}$ are the only edges of $P$ that are adjacent to an unbounded cell of $\ll'$, then $\ll'$ is $2$--critical, and so $\ll$ is $2$--critical (remember, the criticality of $\ll'$ defines the criticality of both $\ll'$ and $\ll$).}
  \label{fig:basic}
  \end{figure}

Let  $Q$ and $p$ be,  respectively, an $n$--gonal region  and a pseudoline of a simple Euclidean arrangement of pseudolines $\ll$.  If  $p$ is adjacent to $Q$, then we say that $p$ {\it is in} $Q$ or that $Q$ {\it is in} $p$.
 If $p$ is in $Q$,  then  we denote the edge of $Q$ in $p$ as $e_{p,Q}$. If $Q$ is formed by $p_1,\ldots ,p_n$, then  the subarrangement of $\ll$ whose pseudolines are $p_1,\ldots ,p_n$ is  the subarrangement of $\ll$ 
{\it  induced} by $Q$. See Figure~\ref{fig:basic}.

If  $p$ and $q$  are pseudolines in an Euclidean arrangement of pseudolines,  then the crossing between $p$ and $q$  will be denoted by $v_{p,q}$. 

In the interest of clarity, and without loss of generality,  we  assume frequently in our arguments that some pseudolines of Euclidean arrangement under consideration are directed. If $p$ is a directed pseudoline  of  an Euclidean arrangement, we denote by $p^+$ (respectively, $p^-$) the semiplane to the right of $p$ (respectively, to the left). Note that if $p$ is a directed pseudoline, $\{p, p^+,p^-\}$ is a partition of $\be$. 

Let $Q$ be a  polygon of an Euclidean arrangement of pseudolines $\ll$ and let  $e$ be an edge of $Q$.
 We say that $e$ is a {\em  critical edge} of  $Q$ if $e$ is adjacent to an unbounded cell of $\ll$. If 
$\ll$ has only one  $(\ge 5)$--gon $P$, and   $k$ is the number of edges of $P$ that are adjacent to an unbounded cell of the subarrangement of $\ll$ induced by the pseudolines in  $P$,  then  we say that $\ll$ is  $k$--{\em critical}. Note that the criticality of $\ll$, is defined by the criticality of the subarrangement of $\ll$ induced by its $(\ge 5)$--gon. See Figure~\ref{fig:basic}.

\newpage

\section{Triangles in arrangements of $\Im$: proof
of Theorem~\ref{igual-t}}\label{somelemas} 

Before proceeding with the proof of Theorem~\ref{igual-t}, we need to establish some results.

\begin {lemma}\label{triangulo-adyacente} 
Let $\ll$ be a simple  Euclidean arrangement of pseudolines. Let  $p,q$ and $r$ be distinct pseudolines of $\ll$ and let   $T$ be the triangular region formed by $p,q$ and $r$. If none of the  pseudolines of $\ll$ crosses $e_{r,T}$, then $\ell\in\{p,q\}$ is in a triangle of $\ll$ and such a triangle is contained in $T$.
\end {lemma}

\noindent{\it Proof.} We prove the theorem with $\ell=p$ (the case $\ell=q$ is  analogous). For brevity of notation,  let $r':=e_{r,T},  p':=e_{p,T}$ and $q':=e_{q,T}$. We proceed by  induction on the number of pseudolines of $\ll$. As $p,q$ and $r$ are distinct pseudolines of $\ll$, it follows that  $|\ll|\ge 3$. It is readily checked that the statement holds for the unique (up to isomorphism) simple Euclidean arrangement with $3$ pseudolines (in this case, $T$ is the required triangle). Thus for some integer $n\ge 3$ we assume (a) the statement holds for every simple Euclidean arrangement of $k\le n$ pseudolines, and (b) that  $|\ll|=n+1.$ 

{\sc Case 1}. There exists a pseudoline $s\in \ll\setminus \{p,q,r\}$ such that $s$ does not cross $T$. Since $s\notin \{p,q,r\}$, then $p,q$ and $r$ are in $\ll':=\ll\setminus \{s\}$. As $\ll'$ is a subarrangement of $\ll$, it follows that $\ll'$ is simple and none of the pseudolines of  $\ll'$ crosses $r'$. Thus, by the inductive hypothesis, there exists a triangle of  $\ll'$, say $T'$, such that $p$ is in $T'$ and $T'$ is contained in  $T$.  As $s$ does not cross $T$ and $T'\subseteq T$, then  $T'$ is a triangle of $\ll$. Moreover, $T'$ is the required triangle.  
 
{\sc Case 2}. Every pseudoline of $\ll \setminus \{p,q,r\}$ crosses $T$. Since none of the pseudolines of $\ll$ crosses $r'$, then every  pseudoline of $\ll\setminus \{p,q,r\}$ crosses $T$ through $p'$ and $q'$.  Note that
 $|\ll|\ge 4$ implies  $\ll\setminus \{p,q,r\}\not=\emptyset$. Let us assume (without loss of generality) that $p$ is directed from $v_{p,q}$ to $v_{p,r}$. Now, we label the pseudolines of  $\ll\setminus \{p,q,r\}$ with $p_1, \ldots ,p_{n-2}$ according to the order in which they are crossed  by $p$. This labeling  is well defined because 
 $\ll$ is simple.

It follows from $p_1\notin \{p,q,r\}$ that  $p,q$ and $r$ are in  $\ll^*:=\ll\setminus \{p_1\}$. Since $\ll^*$ is  
a subarrangement of $\ll$, then $\ll^*$ is simple and none of  the  pseudolines of  $\ll^*$ crosses $r'$. Thus, by inductive hypothesis, there exists a triangle of $\ll^*$, say $T^*$, such that $p$ is in $T^*$  and $T^*$ is contained in $T$.   

Note that if  $\ll=\{p,q,r,p_1\}$, then $T^*=T$ and the triangle formed by $p,q$, and $p_1$ is the required polygon. We may therefore assume that $|\ll|\ge 5$. We also observe that if $p_1$ does not cross $T^*$,  then $T^*$ is the required  triangle of $\ll$. So  we  assume that $p_1$ crosses $T^*$.  

{\sc Subcase 2.1}. $q$ is in $T^*$. As $p, q$ are in $T^*$, and $p_1\notin \ll^*$ and $\ll^*$ is simple,  then $p_2$ is the  third pseudoline in $T^*$. Now,  observe that if   
$v_{p_1,p_2}$ is in $T$, then the triangle formed by $p, p_1$ and $p_2$ is the required  polygon, and in the other case, the triangle formed by $p,q$, and $p_1$ is the required  polygon.
  
{\sc Subcase 2.2}.  $q$ is not in $T^*$.  Since none of the pseudolines of 
$\ll\setminus \{p, q, r\}$ crosses $r'$, then  $r$ is not in $T^*$. Thus, there are distinct integers $i,j\in \{2,\ldots ,n-2\}$ such that  $T^*$ is formed by $p, p_i$ and $p_j$. Without loss of generality, we assume that $i<j$. Since $T^*$ is a triangle of $\ll^*$ and  $p_1$ crosses $T^*$ and $\ll$ is simple,  then when we add $p_1$ to $\ll^*$ to obtain $\ll$, a triangle and a quadrilateral are generated. Both are contained in $T^*$ and  are polygons of $\ll$. Let $Q$ be such a quadrilateral.  As $v_{p,p_1}$ is not in the closed arc of $p$ with endpoints $v_{p, p_i}$ and $v_{p,p_j}$, it follows that $p_1$ crosses $T^*$ through $e_{p_i,T^*}$ and $e_{p_j,T^*}$. Thus $Q$ is formed by $p,p_i, p_j$ and $p_1$. Let $p'_i:=e_{p_i,Q}$, and let $\ll_i$ be the subarrangement of $\ll$ with pseudolines $p,p_1,\ldots ,p_i$ and let $T^i$ be the triangular region formed by $p,p_1$ and $p_i$.  Observe that  $T^i$ is contained in $T$.  Because  $Q$ is a polygon of $\ll$, none of the pseudolines of $\ll$ crosses $p'_i$. In particular none of the pseudolines of $\ll_i$ crosses  $p'_i$. By the induction hypothesis,  $\ll_i$ has a triangle $T^{**}$, such that $p$ is in $T^{**}$ and $T^{**}$ is contained in $T^i$.  Since none of the pseudolines of $\ll$ crosses $p'_i$, and $p_2,p_3,\ldots ,p_{i-1}$ are all the pseudolines of 
$\ll$ that cross $p$ between $v_{p,p_1}$ and $v_{p,p_i}$, it follows that none of the pseudolines of $\ll\setminus \{p , p_1, p_2, \ldots ,p_i\}$ crosses  $T^i$. Thus, $T^{**}$ is a triangle of $\ll$ such that $p$ is in $T^{**}$ and $T^{**}$ is contained in $T^i$  (and hence in $T$), as desired.   
$\square$

\begin {observation}\label{criticidad} 
If $\ll$ is a simple Euclidean arrangement of pseudolines and $P$ is a $(\ge 4)$--gon of $\ll$, then $P$ has at most $2$ critical edges.  
\end {observation}

\begin{figure}[ht]
  \begin{center}
  \includegraphics[width= 4.5 cm, height=3.0cm]{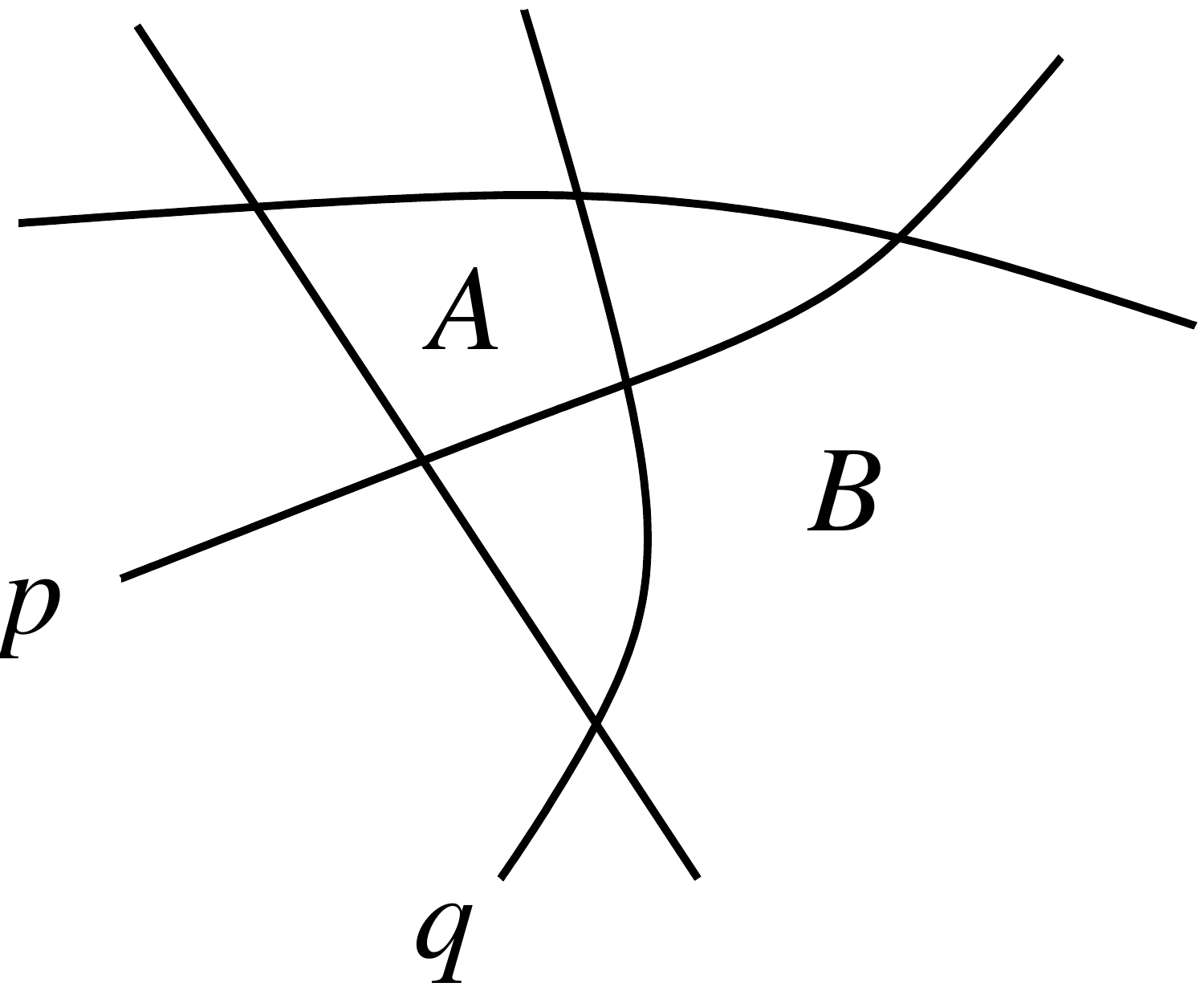}
 \label{unico-cuatro}
  \end{center}
  \caption{\small The unique (up to isomorphism) simple Euclidean arrangement with $4$ pseudolines.}
   \label{fig:arreglo4}
  \end{figure}

\noindent{\it Proof.} We  prove the observation by contradiction. Suppose that $P$ has at least 
$3$ critical edges. Let $\ll'$ be a subarrangement of $\ll$ induced  by $4$ of the pseudolines in $P$, of which at least   $3$  contain critical edges of $P$. Figure~\ref{fig:arreglo4} shows $\ll'$,  the unique (up to isomorphism) simple Euclidean arrangement with $4$ pseudolines. 

\newpage

Since  $\ll'$ is a subarrangement of $\ll$, then $P$ is contained in one  of the eleven  regions  defined by   $\ll'$ (see Figure~\ref{fig:arreglo4}). Since every pseudoline of  $\ll'$ is in $P$, and  $A,  B$ are the only  regions of $\ll'$ that are adjacent to every  pseudoline of $\ll'$, then $P$ must be contained in one of  $A$ or $B$. Observe that in both cases, neither $p$ nor $q$ can contain a critical edge  of $P$. This contradicts our assumption that $\ll'$ has at least $3$ critical edges of $P$.  $\square$

\begin {proposition}\label{almenos-t} 
If  $\ll\in \Im$  then  every non-critical edge  of  the  $(\ge 5)$--gon of  $\ll$ is also an edge of  a  triangle of  
$\ll$.
\end {proposition}
 
\noindent{\it Proof.} Let $P$ be the  $(\ge 5)$--gon of $\ll$. Let $e'_1$ be any non-critical edge of $P$ and let $e'_2, e'_3$ be the edges of  $P$ that are adjacent to $e'_1$. For 
$i=1,2,3$ let  $e_i$ be the  pseudoline of $\ll$ that contains $e'_i$. Without loss of generality we may assume that each $e_i$ is directed in such a way that $P$ lies to the right of $e_i$. 

Since $e'_1$ is a non-critical edge of $P$, then the region to the left of $e'_1$ is a polygon of  $\ll$. We denote by $Q$ such a polygon. $Q$ must be a triangle or a quadrilateral, because $P$ is the only $(\ge 5)$--gon of $\ll$. Assume, for the purpose of contradiction, that  $Q$ is a quadrilateral. 

Before proceeding,  note that $e_j$  ($j=2,3$) contains an edge of $Q$ because  $\ll$ is simple.  Thus, as $Q$ is quadrilateral, $\ll\setminus \{e_1, e_2, e_3\}$ contains a  pseudoline $x$, such that $Q$ is formed by $x, e_1, e_2$ and $e_3$. Let $x':=e_{x,Q}$ and let  $x_2$  and $x_3$ be the two unbounded subarcs of $x$ obtained by deleting $x'$. Without loss of generality, we assume that  $x_j$ ($j=2,3$) is the unbounded subarc of $x$ that has endpoint at $v_{x,e_j}$. Since $v_{x,e_j}$ is the unique intersection point between  $x$ and $e_j$, then   $x_j$  is contained in $e_j^-$. Thus,  $x_j $ is not in $P$. On the other hand, is clear that $P$  cannot be adjacent to $x'$, and hence $P$ is not in $x$, which contradicts the fact that every pseudoline of $\ll$ is in $P$.  $\square$  

\begin {proposition}\label{alomas-t} 
If  $\ll \in \Im$ then every triangle of  $\ll$ has a common edge with the $(\ge 5)$--gon of $\ll$.
\end {proposition}
 
\begin{figure}[ht]
  \begin{center}
  \includegraphics[width= 3.5 cm, height=2.5cm]{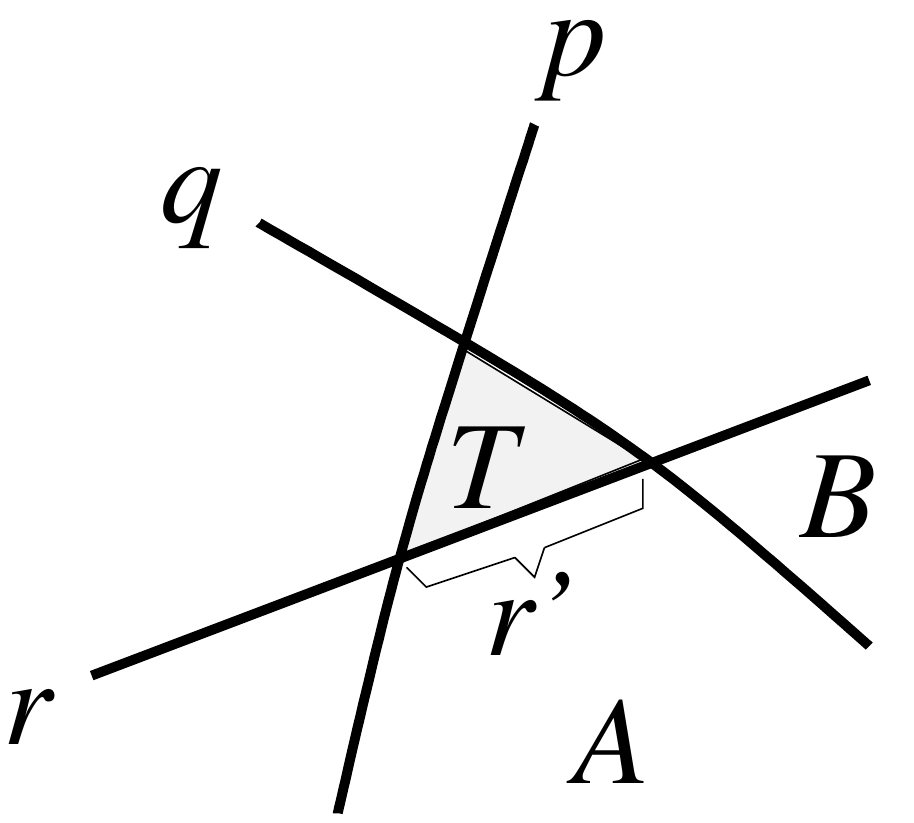}
  \end{center}
  \caption{\small The unique (up to isomorphism) simple Euclidean arrangement with $3$ pseudolines.}
  \label{fig:arreglo3}

  \end{figure}

\noindent{\em Proof.} Let $T$ be any  triangle of $\ll$ and let  $p,q$ and $r$ be the pseudolines in $T$. We denote by  $P$ the $(\ge 5)$--gon of $\ll$. Let $\ll'$ be the subarrangement of $\ll$ induced by $T$. Figure~\ref{fig:arreglo3} shows $\ll'$,  the unique  (up to isomorphism) simple Euclidean arrangement with $3$ pseudolines.

Since  $\ll'$ is a subarrangement of  $\ll$, then $P$ must be contained in one of the  six unbounded regions   defined by  $\ll'$ (see Figure~\ref{fig:arreglo3}). By symmetry, we may assume  that $P$ is contained in $A$ or  $B$. Since $\ll \in \Im$, $p$ is in $P$ and  hence $P$ cannot be contained in $B$. So $P\subset A$. Let $r':=e_{r,T}$. Because $P\subset A$ and $r$ is in  $P$,  there exists a subarc $r''\subseteq r'$ that is an edge of $P$. On the other hand, none of the pseudolines of  $\ll$ crosses  $r'$ because  $T$ is a polygon of $\ll$. Thus,  $r''=r'$, as desired. $\square$
  
\begin {theorem}\label{igual-t} 
If $\ll\in \Im$  and $\ll$ is $k$--critical,  then  the number of triangles of $\ll$ is $|\ll|-k$.
\end {theorem}

\noindent{\it Proof.} Let $P$ be the $(\ge 5)$--gon of $\ll$ and $n:=|\ll|$. Since the  pseudolines of
$\ll$ intersect  pairwise exactly once, it follows that distinct edges of $P$ come from distinct pseudolines of $\ll$. This implies that $P$ has at most  $n$ edges. On the other hand, since every pseudoline of $\ll$  is in $P$,  then $P$ has at least  $n$ edges. Thus, $P$ is an $n$--gon and the assertion follows from Propositions~\ref{almenos-t} and~\ref{alomas-t}. $\square$

\section{Triangles in simple Euclidean arrangements with one $(\ge 5)$--gon: proof
of Theorem~\ref{igual-ss}}\label{general} 

We begin the proof of Theorem~\ref{igual-ss} with the 
Proposition~\ref{coro-triangulo-adyacente}.

\begin {proposition}\label{coro-triangulo-adyacente} 
Let $\ll$ be a simple Euclidean arrangement of $n\ge 3$ pseudolines. Then every pseudoline of
$\ll$ is adjacent to at least one triangle.  
\end {proposition}

\noindent{\it Proof.} We apply induction on $|\ll|$. The assertion is  trivial for the unique (up to isomorphism) simple Euclidean arrangement with $3$ pseudolines. Thus for some integer $n\ge 3$ we assume (a) the statement holds for every simple Euclidean arrangement of $k\le n$ pseudolines, and (b) that  $|\ll|=n+1.$

Let $t$ be any pseudoline of  $\ll$ and let $s$ be a pseudoline of  $\ll\setminus \{t\}$.
 By the inductive hypothesis, $\ll':=\ll\setminus \{s\}$  has a triangle $T$, such that $t$ is in $T$.  If $s$ does not cross $T$,  then   $T$ is a triangle of $\ll$ and   $T$ is the required triangle.  So we assume that $s$ crosses  $T$.  Let $x,y\in \ll'$ be the other two pseudolines in $T$.
Observe that if $s$ crosses  $T$ through  $t$ and  $p\in \{x,y\}$, then the triangle formed by $t, s$ and $p$ is the required triangle.   We may therefore assume that $s$ crosses $T$ through  $x$ and $y$. This implies that none of the pseudolines of $\ll$ crosses the subarc of  $x$ with endpoints $v_{x,t}$ and $v_{x,s}$. Hence, by Lemma~\ref{triangulo-adyacente}, there exists a triangle of $\ll$, say $T_t$,  such that 1) $T_t$  is contained in the triangular region formed by  $x,s$ and $t$, and 2) $t$ is in $T_t$,  as desired.  $\square$

\begin{figure}[ht]
  \begin{center}
  \includegraphics[width= 5 cm, height=3.0cm]{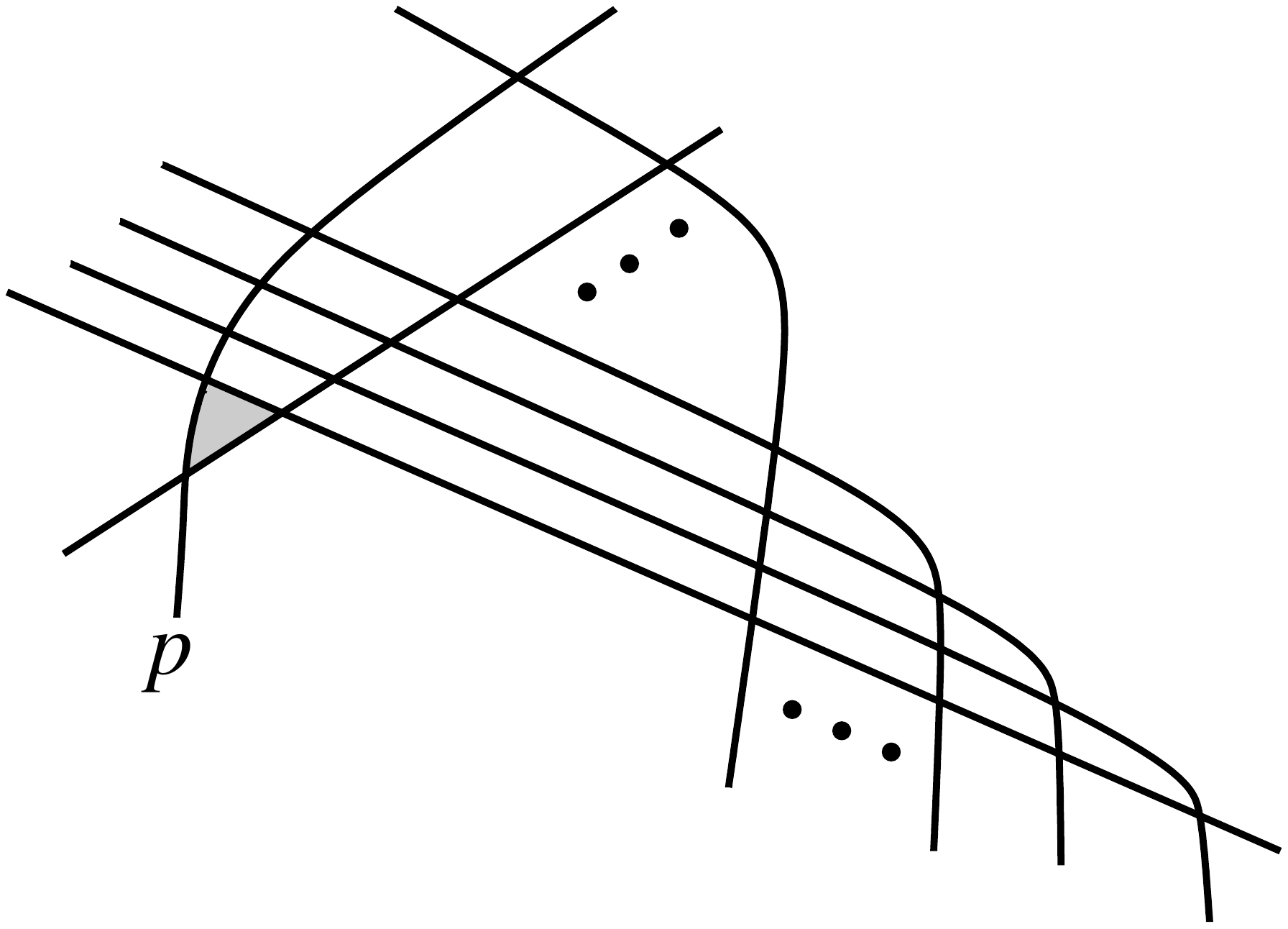}
  \end{center}
  \caption{\small Pseudoline $p$ is adjacent to  exactly one triangle.}
  \label{fig:ccoro}
  \end{figure}

The construction in Figure~\ref{fig:ccoro} shows that Proposition~\ref{coro-triangulo-adyacente} is best possible.

\begin {lemma}\label{n-gono-adyacente} 
Let $\ll$ be a simple Euclidean arrangement of  pseudolines and let $\mm$ be a subarrangement of $\ll$.  If $Q$ is a  $(\ge5)$--gon of  $\mm$ and $w$ is an edge of  $Q$ such that none of the pseudolines of   $\ll$ crosses $w$, then every edge of $Q$ that is adjacent to $w$ contains a subarc that is an edge of a  $(\ge5)$--gon of $\ll$. Moreover, $Q$ contains such a $(\ge5)$--gon.   
\end {lemma}
 
\noindent{\it Proof.} Let $q$ be any of the two edges of $Q$ that are adjacent to 
$w$.  Let us assume  (without loss of generality) that $q$ is directed in such a way that $Q$ lies to the right of $q$.  Let $p\not=w$ (respectively $r\not=q$) the edge of $Q$ that is adjacent to $q$ (respectively $w$). Since  $Q$ is a $(\ge 5)$--gon, $p$ is not adjacent to $r$. 

Let $\cc$ be the set of all pseudolines of $\ll$ that cross $q$, and let $\cc'$ be the set of all pseudolines of  $\ll\setminus ( \mm \cup \cc)$. Let $\mm'$ be the subarrangement of  $\ll$  obtained by deleting the pseudolines in  $\cc$. 
 Let  $Q'$ be the polygon of  $\mm'$ that is adjacent to $q$ and lies to the right of $q$. Clearly, $Q'$ is contained in $Q$.  Note that both  $q$ and  $w$ are edges of  $Q'$. This implies that  $v_{q,p}, v_{q,w}$ and $v_{w,r}$ are vertices of $Q'$. Since $\ll$ is simple and $q$ (respectively $w$) is an edge of $Q'$, then there exists a subarc $p'\subseteq p$ (respectively $r'\subseteq r$) with an endpoint at $v_{q,p}$ (respectively  $v_{w,r}$) such that $p'$ (respectively $r'$) is an edge of $Q'$. Since  $p$ is not adjacent to $r$,  $p'$ is not adjacent to $r'$. We denote by  $v_{r'}$ (respectively $v_{p'}$) the vertex of  $Q'$ that is an endpoint of  $r'$ (respectively $p'$) but not an endpoint of $w$ (respectively $q$).
 
 Since the edges of  $Q'$ form a cycle, then there is a path $\Gamma$, formed by edges of $Q'$ with endpoints  $v_{p'}$ and  $v_{r'}$,   such that $\Gamma\cap \{w,q\}=\emptyset$. As $p'$ and $r'$ are not adjacent, then  $\Gamma$ contains at least one edge, and hence $Q'$ is a $(\ge 5)$--gon. 
 
If $|\cc|=0$ then  $Q'$ is a  $(\ge 5)$--gon of $\ll$, and $q$ and $Q'$ are as desired. So we assume that 
 $|\cc|=n>0$.

We label the pseudolines of  $\cc$ with $q_1, \ldots ,q_n$ according to the order in which they are found when $q$ is walked from $v_{q,w}$ to $v_{q,p}$. This labeling  is well defined because  $\ll$ is simple.

For $i=1, \ldots ,n-1$ let $\alpha_i$ be the subarc of  $q$ with endpoints $v_{q,q_i}$ and $v_{q,q_{i+1}}$. Let  $\alpha_0$ be the subarc of  $q$ with endpoints $v_{q,w}$ and  $v_{q,q_{1}}$ and let $\alpha_n$ be the  subarc of $q$ with endpoints $v_{q,q_{n}}$ and  $v_{q,p}$. For $j=0,\ldots ,n$  (without loss of generality we assume that  $\alpha_j$ has the same direction as  $q$)  let $Q'_j$ be the polygon  to the right of  $\alpha_j$ ($\alpha_j$ is an edge of $Q'_j$). Clearly, $Q'_j$ is contained in $Q'$. From the definition of $Q'_j$ it is not difficult to see that $Q'_j$ is a polygon of $\ll$. Thus, if $Q'_j$ is a $(\ge 5)$--gon for some $j\in \{0,\ldots ,n-1\}$ we are done. So assume that every element of $\{Q'_0, \ldots ,Q'_{n-1}\}$ is a $(\le 4)$--gon. 
 
On the other hand, observe that for   $i=1,\ldots ,n-1$, $q_i$ and  $q_{i+1}$ are in $Q'_i$. Moreover, note  that  $e_{q_i,Q'_{i-1}}=e_{q_i,Q'_{i}}$. For brevity of notation let $q'_i:=e_{q_i,Q'_{i-1}}=e_{q_i,Q'_{i}}$.  

Since none of the pseudolines of  $\cc$ crosses $w$ and $\ll$ is simple, then $w, \alpha_0, q'_1$  and  a subarc $r_0$ of $r'$, which has  a common endpoint  with  $w$,  belong to the  edge set  of $Q'_0$. Thus,  $Q'_0$ is a quadrilateral with  edge set  $\{ w, \alpha_0, q'_1, r_0\}$. Since $\ll$ is simple,  $r_0$ is a proper subarc of  $r'$ and consequently  $q_1$ must cross $r'$.  Let  $v_1$ be the crossing  between $r'$ and  $q_1$ ($v_1$ is common endpoint of both $r_0$ and $q'_1$). On the other hand, since  $Q'_0$ is a quadrilateral of $\ll$, then none of the pseudolines of  $\cc$ crosses $q'_1$. 

The argument in the next paragraph works for $i=1,\ldots ,n-1$ if we apply it repeatedly in ascending order.

As $\ll$ is simple and none of the pseudolines of $\cc$ crosses $q'_i$, then $q'_i, \alpha_i, q'_{i+1}$ and a subarc $r_i$ of $r'\setminus (r_0\cup \ldots \cup r_{i-1})$, which has an endpoint at $v_i$,  belong to the  edge set  of $Q'_i$. Thus,  $Q'_i$ is a  quadrilateral with edge set  $\{q'_i, \alpha_i, q'_{i+1},r_i\}$. As $\ll$ is simple,  then $r_i$ is a proper subarc of  $r'\setminus  (r_0\cup \ldots \cup r_{i-1})$ and consequently  $q_{i+1}$ must cross $r'\setminus  (r_0\cup \ldots \cup r_{i-1})$. Let  $v_{i+1}$ be the intersection point between $r'\setminus  (r_0\cup \ldots \cup r_{i-1})$ and  $q_{i+1}$ ($v_{i+1}$ is common endpoint of both $r_i$ and  $q'_{i+1}$). On the other hand, as $Q'_i$ is a quadrilateral of $\ll$, then none of the pseudolines of  $\cc$ crosses $q'_{i+1}$.
 
 From the previous   paragraph we know that $Q'_{n-1}$ is a quadrilateral with edges  $q'_{n-1}, \alpha_{n-1}, q'_n$ and $r_{n-1}$.  Since $\ll$ is simple and $r_{n-1}$ is a proper subarc of  $r'\setminus  (r_0 \cup \ldots \cup r_{n-2})$, which has an endpoint at $v_{n}$, then the subarc $r_n:=r'\setminus  (r_0\cup \ldots \cup r_{n-1})$ whose endpoints are $v_n$ and  $v_{r'}$, is an edge of  $Q'_n$. Thus,  $q'_n, \alpha_n, p'$,  all the edges of $\Gamma$, and  $r_n$ are the edges of   $Q'_n$. As $\Gamma$ has at least one edge and $\Gamma\cap \{q'_n, \alpha_n, p', r_n\}=\emptyset$, then $Q'_n$ is a $(\ge 5)$--gon, as desired.   $\square$

\vskip 0.2cm

\noindent{\it Proof of Theorem~\ref{igual-ss}}. We proceed by induction on the number of pseudolines of $\ll$ that are not in $P$. If the number of pseudolines of $\ll$ that are not in $P$ is zero, then $\ll\in \Im$ and we are done by Theorem~\ref{igual-t}.  So we assume that  $\ll$ has $\eta \ge 1$ pseudolines that are not in  $P$ and that the assertion holds for  every simple Euclidean arrangement $\nn$  with one $(\ge 5)$--gon  such that $\nn$ has $j<\eta$  pseudolines that are not in the  $(\ge 5)$--gon of $\nn$. 

Let $p$ be a pseudoline of $\ll$ that is not in $P$ and let $\ll'$ be the subarrangement of $\ll$ obtained by deleting $p$. It is not difficult to see that  $\ll'$ is a simple Euclidean arrangement with one $(\ge 5)$--gon.  In fact, note that $P$ is the $(\ge 5)$--gon of $\ll'$. Then both  $\ll$ and $\ll'$ are $k$--critical and by the inductive hypothesis $\ll'$ has $|\ll'|-k=(|\ll|-1)-k=n-1-k$ triangles.

Now we partition the set of triangles of  $\ll$ into two subsets. The first subset is denoted by $A$ and consists of all the triangles that have an edge in $p$. The second subset is denoted by $B$ and contains the rest of the triangles of  $\ll$. So the number of triangles of $\ll$ is $|A|+|B|$.

It follows from the definition of $\ll'$  that every triangle of  $B$ is a triangle of $\ll'$. On the other hand, it is clear that if $Q'_1,\ldots ,Q'_{\ell}$ are the polygons of  $\ll'$ that are crossed by $p$, and $\alpha$ of them are triangles, then the number of triangles of $\ll'$ is $|B|+\alpha$. To complete the proof, it suffices to show that   $|A|-1=\alpha$.

Let $A_1$ be the subset of  $A$ consisting  of all the triangles that are contained in a polygon of 
$\ll'$, and let $A_2:=A\setminus A_1$. Note that  $A_2$ is the set of all  triangles of $\ll$ that have a critical edge in $p$. Without loss of generality we may assume that $p$ is directed.  For $j\in\{1,\ldots ,\ell\}$ we denote by $p_j$ the open subarc of $p$ that intersects $Q'_{j}$ and by 
$Q_{l,j}$ (respectively,  $Q_{r,j}$)  the polygon of $\ll$ to the left  (respectively, right) of $p_j$ such that  $p_j=e_{p,Q_{l,j}}=e_{p,Q_{r,j}}$. Observe that  $\{Q_{l,j},  Q_{r,j}, p_j\}$ is a partition of $Q'_j$. From the  definitions of  $A$,  $Q_{l,j}$, $Q_{r,j}$ and $A_1$, it is not difficult to see that $A_1=A\cap \{Q_{l,1}, Q_{r,1},\ldots ,Q_{l,\ell},  Q_{r,\ell}\}$.  Because $P$ is the unique $(\ge 5)$--gon of $\ll$ and $p$ is not in $P$, every $Q_{d,j}$ is a triangle or a quadrilateral for $d\in \{l,r\}$. Moreover, since $\ll$ is simple and  the  pseudolines of $\ll$ intersect  pairwise exactly once,  $Q_{l,j}$ and $Q_{r,j}$ are not  both triangles. Thus, $Q'_j$ is a triangle if and only if exactly one of  $Q_{l,j}$ or $Q_{r,j}$ is a triangle.  This implies that $|A_1|=\alpha$. Thus, all we need to  show is that $|A_2|=1$.

We divide the remainder of the proof into two parts.

\begin{figure}[ht]
  \begin{center}
  \includegraphics[width= 8cm, height=4cm]{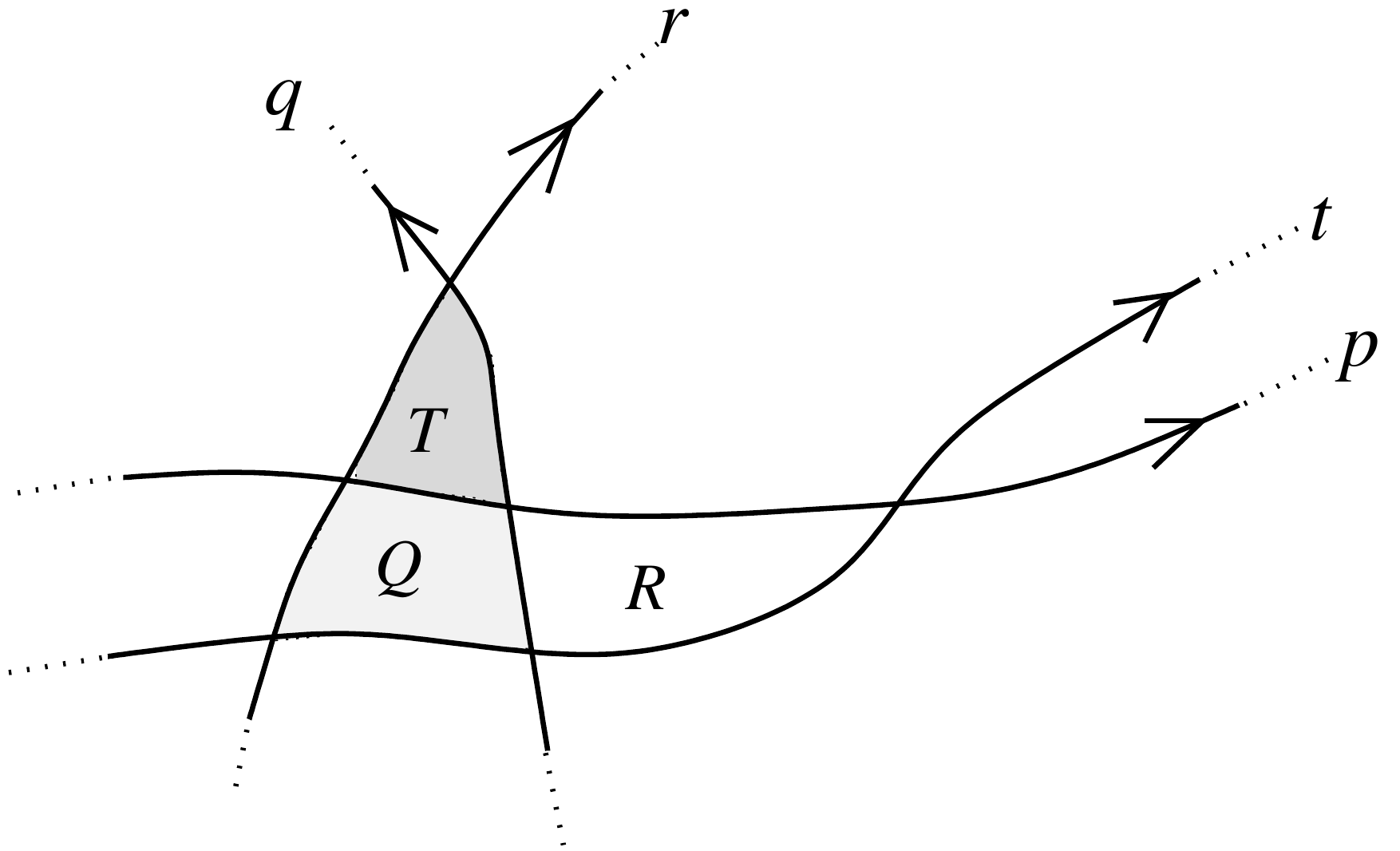}
  \end{center}
  \caption{\small The unique (up to isomorphism) simple Euclidean arrangement with $4$ pseudolines.}
  \label{solocuatro}
  \end{figure}

{\bf Part I}. In this part we shall show that  $|A_2|\ge 1$. To obtain a contradiction (that $P$ is in $p$) we assume that $A_2$ is empty. By Proposition~\ref{coro-triangulo-adyacente}, $p$ is in at least one triangle, say $T$,  of  $\ll$. $T$ must be an element of $A_1$ because $A_2$ is empty. Thus, there exists an $m\in\{1,\ldots ,\ell\}$ such that $T=Q_{l,m}$ or $T=Q_{r,m}$. Without loss of generality  (changing the direction of $p$ if necessary)  we may assume that  
$T:=Q_{l,m}$. So $Q:=Q_{r,m}$ is a quadrilateral.  Figure~\ref{solocuatro}  shows  the unique (up to isomorphism) simple Euclidean arrangement with $4$ pseudolines. So we may assume (without loss of generality)  that $p$ and the three pseudolines in $Q'_{m}$ are directed and labeled as in Figure~\ref{solocuatro}.

{\sc Case 1.} $R$ is a triangle of  $\ll$. This implies that none of the pseudolines of $\ll$ crosses $e_{p,R}$. Analogously, since $T$ is a triangle of  $\ll$, none of the pseudolines of $\ll$ crosses $e_{q,T}$. Since $A_2$ is empty,  the region $S$ to the left of  $e_{p,R}$ must be a polygon of $\ll$.   Because $\ll$ is simple,  $v_{t,p}, v_{q,p}$ and $v_{q,r}$ are vertices of $S$ (and hence $e_{p,R}$ and $e_{q,T}$ are  edges of $S$). In this circumstances, it is easy to see that $S$ must be a $(\ge 5)$--gon. But  $P$ is the unique $(\ge 5)$--gon of $\ll$. So $P=S$, and therefore $P$ is in $p$.

{\sc Case 2.} $R$ is not a triangle of $\ll$. This means that at least one of the pseudolines of $\ll^*:=\ll\setminus \{p,q,r,t\}$ crosses $R$. Since $Q$ is a polygon of  $\ll$ and $\ll$ is simple,  none of the pseudolines of  $\ll$ crosses $e_{q,Q}\cup \{v_{q,t}, v_{p,q}, v_{t,p}\}$. Then $x\in \ll^*$ crosses $R$  if and only if $x$ crosses $e_{p,R}$ and $e_{t,R}$.  Let  $a\in \ll^*$ be the first pseudoline that crosses $p$ when $p$ is walked from $v_{p,t}$ to $v_{p,q}$  ($a$ is well-defined because $\ll$ is simple).

{\sc Subcase 2.1.} $\ll^*$ contains at least one element, say $b$, such that the intersection point between $a$ and $b$ is  in $R$. In such a case, note that a pentagon $S$, is formed by $a,b,p,q$ and 
$t$. Thus, $S$ is a $(\ge 5)$--gon of the subarrangement of $\ll$ induced by $a,b,p,q$ and $t$, and $e_{q,Q}$ is an edge of $S$. By Lemma~\ref{n-gono-adyacente}, $e_{p,S}$ contains a subarc  that is an edge of a $(\ge 5)$--gon of $\ll$, say $S'$. But  $P$ is the unique $(\ge 5)$--gon of $\ll$. So $P=S'$, and therefore $P$ is in $p$. 

\begin{figure}[ht]
  \begin{center}
  \includegraphics[width= 8cm, height=4cm]{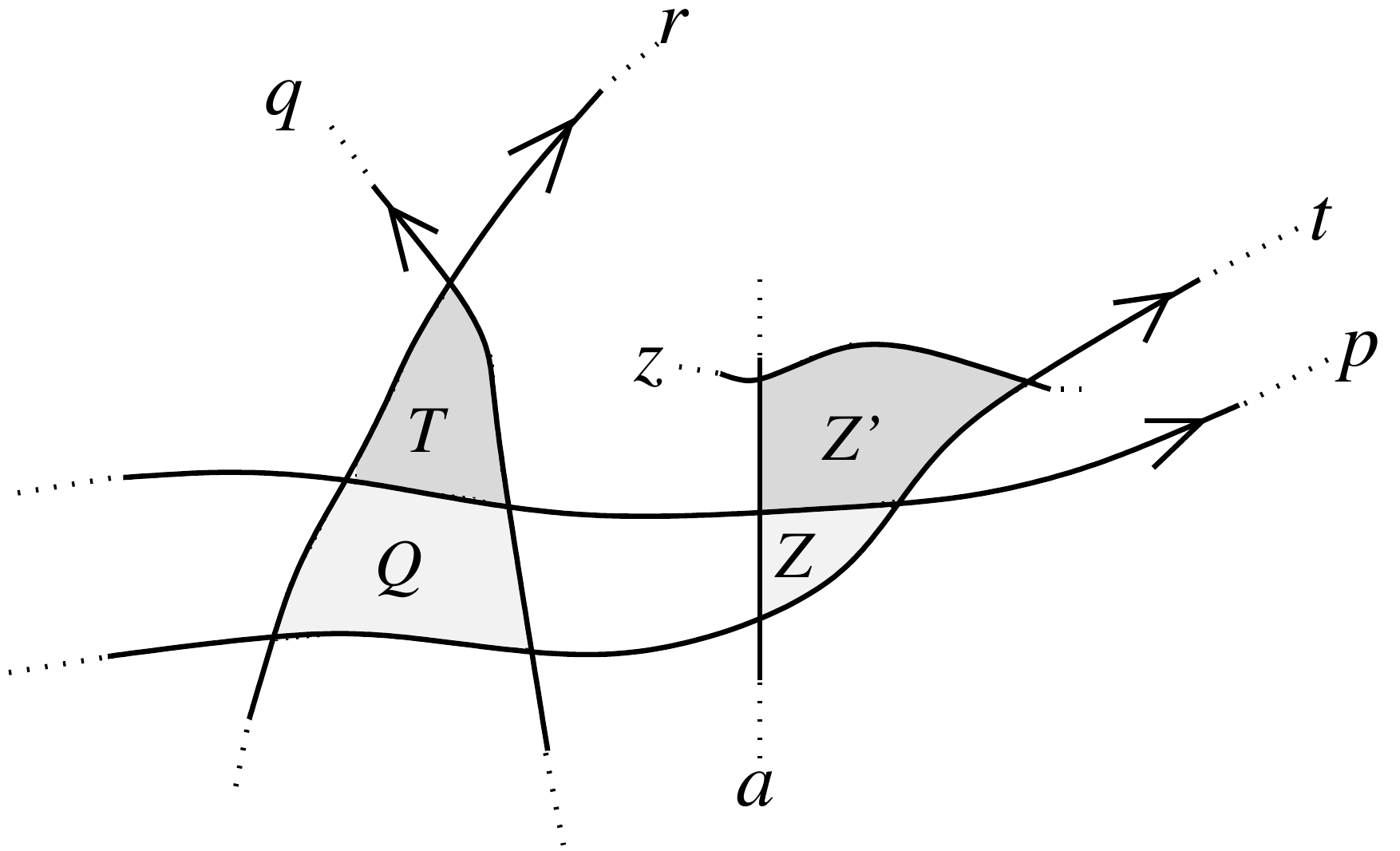}
  \end{center}
  \caption{\small The form of  the subarrangement of 
   $\ll$ induced by  $a,p,q,r,t$ and $z$ if no pseudoline of $\ll^*$  crosses $a$ in $R$.}
  \label{subdos}
  \end{figure}

{\sc Subcase 2.2.} None of the pseudolines of  $\ll^*$ crosses $a$ in $R$. This implies that the triangular region $Z$, formed by  $a,p$, and $t$ is a triangle of $\ll$ (see Figure~\ref{subdos}), and hence none of the pseudolines of $\ll$ crosses $e_{p,Z}$. Since $A_2$ is empty, the region $Z'$ to the left of $e_{p,Z}$ must be a polygon of $\ll$. Then $Z'$  is a triangle or a quadrilateral, because $P$ is the unique $(\ge 5)$--gon of $\ll$ and $p$ is not in $P$. On the other hand, since in a simple Euclidean arrangement there are no triangles with a common edge,  $Z'$ is a quadrilateral. Because $\ll$ is simple, $Z'$ has an edge in each of  $a$ and $t$. Let $z$ be the fourth pseudoline in $Z'$. As $v_{q,t}, v_{r,t}\in p^+$ and  $v_{z,t}\in p^-$, then  $z\notin\{q,r\}$. Observe that if $v_{p,z}\in  e_{p,R}$ then $v_{t,z}\in  e_{t,R}$, because $e_{q,R}=e_{q,Q}$ cannot be crossed and $\ll$ is simple. But $v_{t,z}\in p^-$ and $e_{t,R}\in p^+$,  so $v_{p,z}\notin e_{p,R}$. Moreover, since $\ll$ is simple and  $Q$ is a polygon of $\ll$, if   $p':=e_{p,R} \cup e_{p,Q}\cup  \{v_{p,r}, v_{p,q}, v_{p,t}\}$, then $v_{p,z}\notin p'$.

The last two cases from Part I now arise, depending on whether the intersection point between $r$ and  $z$ occurs in  $p^+$ or  $p^-$.

\begin{figure}[ht]
  \begin{center}
  \includegraphics[width= 15cm, height=5.0cm]{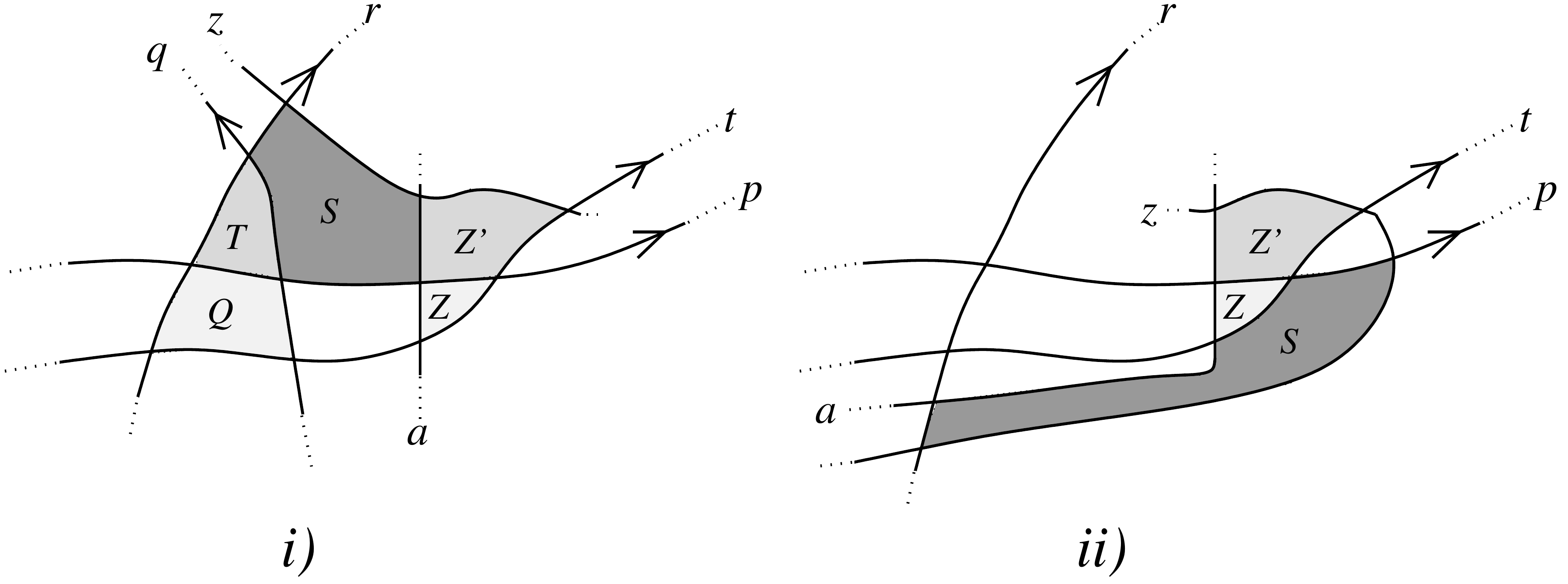}
  \end{center}
  \caption{\small  $i)$ illustrates the case  where the intersection point between $r$ and $z$ is in $p^-$; and $ii)$ illustrates the case  where the intersection point between $r$ and $z$ is in $p^+$.}
\label{subsub}
  \end{figure}

{\sc Subcase 2.2.1} $v_{r,z}\in p^-$.  Because   $\ll$ is  simple and $T$ is a polygon of  $\ll$,    
$v_{z,q}\notin e_{q,T}\cup \{v_{p,q},v_{q,r}\}$. So a pentagon $S$, is formed by $a,p,q,r$ and $z$. 
Moreover,  $S$  is a $(\ge 5)$--gon of the subarrangement of $\ll$ induced by $a,p,q,r$ and $z$, and 
$e_{q,T}=e_{q,S}$,  see  Figure~\ref{subsub} $i)$. Since  $T$ is a polygon of $\ll$, none of the 
pseudolines of $\ll$ crosses $e_{q,S}$. Thus, by  Lemma~\ref{n-gono-adyacente}, $e_{p,S}$ contains a subarc  that is an edge of a $(\ge 5)$--gon of $\ll$, say $S'$. But  $P$ is the unique $(\ge 5)$--gon of $\ll$. So $P=S'$, and therefore $P$ is in $p$.

{\sc Subcase 2.2.2} $v_{r,z}\in p^+$. Then $v_{p,z}\in t^+$ because $v_{p,z}\notin p'$. Together, $v_{r,z}\in p^+$ and  $v_{p,z}\in t^+$ imply that the triangular region $H$ formed by $r,t$ and $z$ is contained in  $t^+$.  Because both $a$ and  $p$ cross $t$ in $e_{t,H}$, $a$ and  $p$ cross $H$. Since $v_{a,z}\in t^-$,  it follows that $a$ crosses  $H$ through $e_{t,H}$ and $e_{r,H}$. Analogously,  $p$ crosses  $H$  through $e_{t,H}$ and $e_{z,H}$ because $v_{p,r}\in t^-$. It follows from the last two assertions and the fact that $v_{a,p}$ is not in $H$ that $a,p,r,t$ and $z$ form a pentagon $S$. 
Moreover,  $S$  is a $(\ge 5)$--gon of the subarrangement of $\ll$ induced by $a,p,r,t$ and $z$, and 
$e_{t,Z}=e_{t,S}$,  see  Figure~\ref{subsub} $ii)$. Since  $Z$ is a triangle of $\ll$, none of the 
pseudolines of $\ll$ crosses $e_{t,S}$. Thus, by  Lemma~\ref{n-gono-adyacente}, $e_{p,S}$ contains a subarc  that is an edge of a $(\ge 5)$--gon of $\ll$, say $S'$. But  $P$ is the unique $(\ge 5)$--gon of $\ll$. So $P=S'$, and therefore $P$ is in $p$.

\begin{figure}[ht]
  \begin{center}
  \includegraphics[width= 11cm, height=7cm]{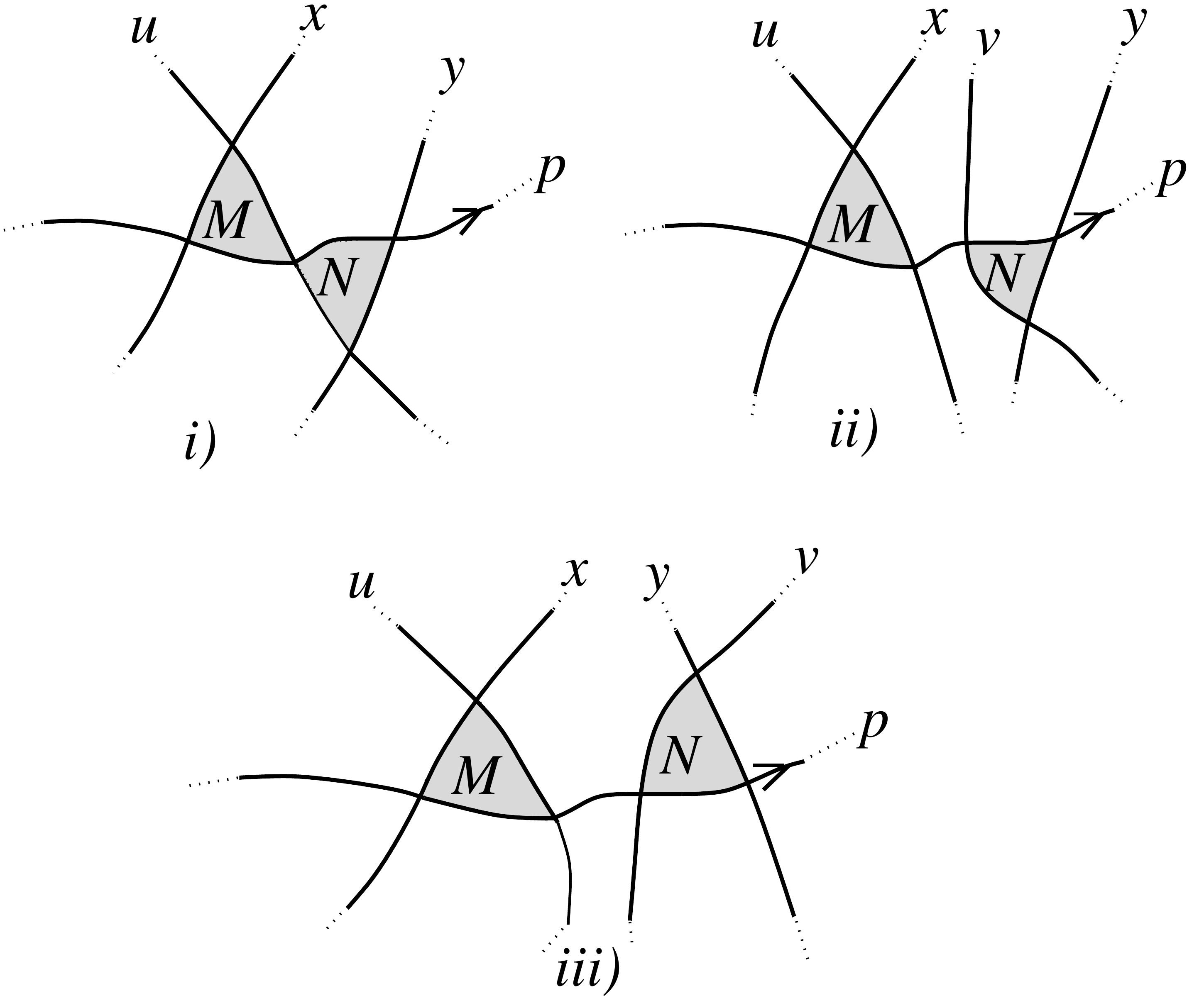}
  \end{center}
  \caption{\small There are only three (up to isomorphism) possible configurations for $M, N$ and 
  $p$.}
  \label{trescasos}
  \end{figure}

\newpage
{\bf Part II}. In this part we shall prove that $|A_2|\le 1$. Suppose not. Let $M, N$ be distinct triangles of   $A_2$. It is not difficult to see that the configurations shown in Figure~\ref{trescasos} are all  (up to isomorphism)  the possible configurations where both $M$ and $N$ are in $p$; note that it is impossible for $M$ and $N$ to be on the same side of $p$ and share a vertex, since $\ll$ is simple. 

\begin{figure}[ht]
  \begin{center}
  \includegraphics[width= 7cm, height=4.0cm]{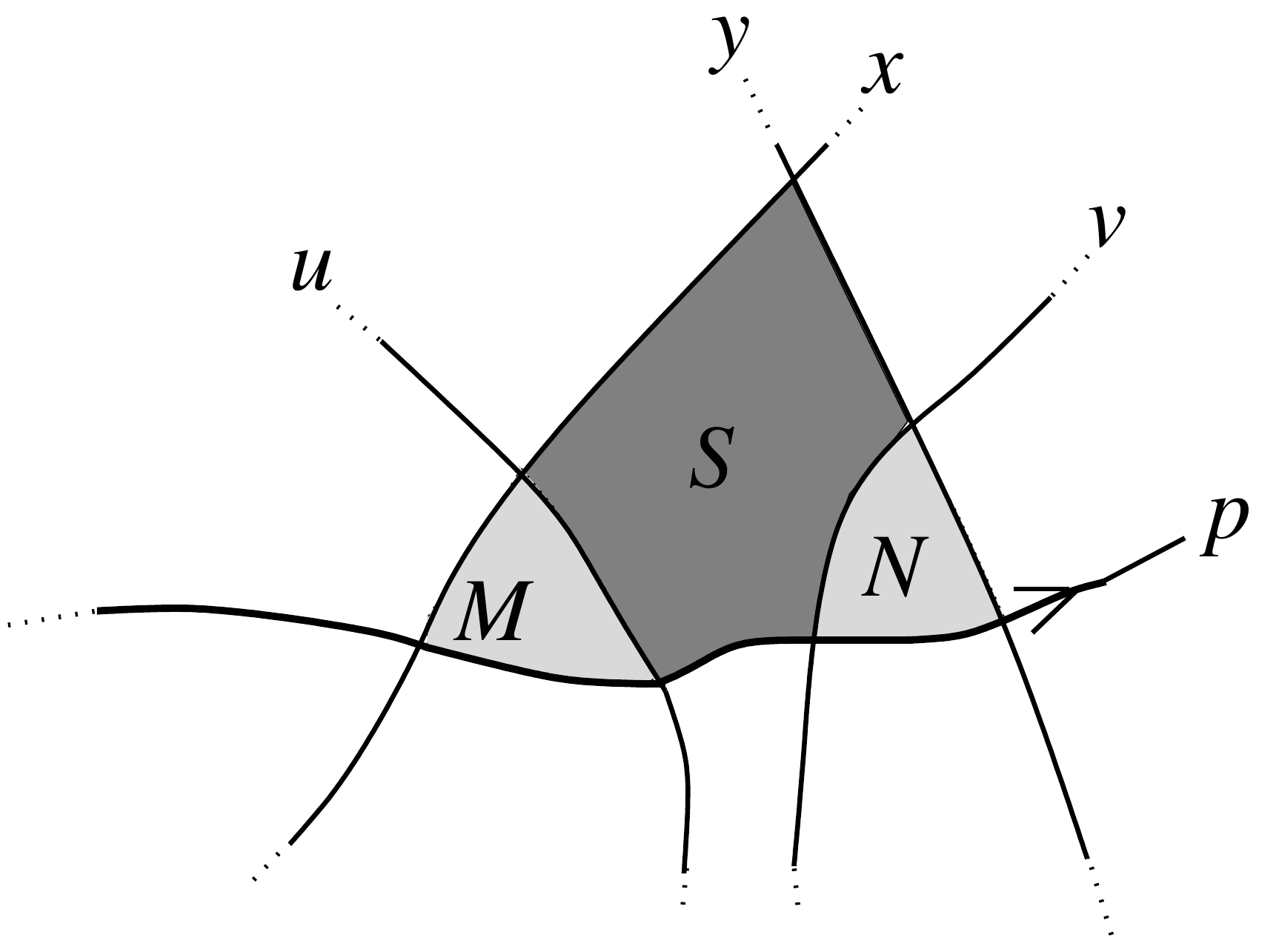}
  \end{center}
  \caption{\small The unique (up to isomorphism)  configuration possible for $M, N$ and 
  $p$.}
  \label{elunico}
  \end{figure}

Assume that  $M, N$ and $p$  are placed as in Figure~\ref{trescasos} $i)$ or $ii)$.   Since $\ll$ is simple, $v_{x,y}\in p^+$ or $v_{x,y}\in p^-$. 
If $v_{x,y}\in p^+$ (respectively, $v_{x,y}\in p^-$) then $e_{p,M}$ (respectively, $e_{p,N}$) is not a critical edge of $M$ (respectively, $N$), since, e.g., $N$ is contained in the triangle  formed by $u$, $x$ and $y$. However, this contradicts the fact that every triangle in $A_2$ has critical edge in $p$.  

So we may assume that $M, N$ and $p$ are placed as in Figure~\ref{trescasos} $iii)$. Because $e_{p,M}$ and $e_{p,N}$ are, respectively, critical edges of   $M$ and $N$, $v_{x,y}\in p^-$. This implies that 
a pentagon, $S$, is formed by  $p,x,y,u,$ and $v$ (see Figure~\ref{elunico}).
   
Moreover,  $S$  is a $(\ge 5)$--gon of the subarrangement of $\ll$ induced by $p,x,y,u$ and $v$, and 
$e_{u,S}=e_{u,M}$. Since  $M$ is a triangle of $\ll$, none of the 
pseudolines of $\ll$ crosses $e_{u,S}$. Thus, by  Lemma~\ref{n-gono-adyacente}, $e_{p,S}$ contains a subarc  that is an edge of a $(\ge 5)$--gon of $\ll$, say $S'$. But  $P$ is the unique $(\ge 5)$--gon of $\ll$. So $P=S'$, and therefore $P$ is in $p$. $\square$

\section{ Arrangements in $\Im$ are stretchable: proof of Theorem~\ref{rectif}}\label{secrec}
Recall that Theorem~\ref{rectif} asserts that if $\ll$ is an arrangement in $\Im$, then $\ll$ is stretchable. 
The proof is by induction on $|\ll|$. In~\cite{goodman-pollack}, Goodman and Pollack showed that any Euclidean arrangement with $8$ pseudolines is stretchable. So we may assume that $|\ll|=n \ge 9$ 
 and that every arrangement of  $\Im$ with $j<n$ pseudolines is stretchable. 

\begin{figure}[ht]
  \begin{center}
  \includegraphics[width= 15cm, height=3.5cm]{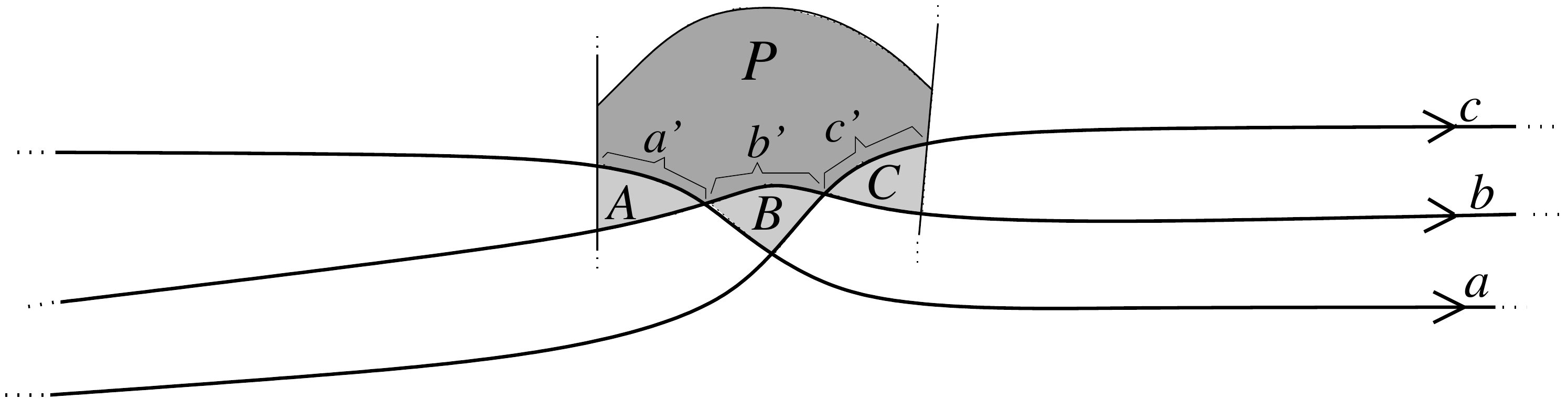}
  \end{center}
\caption{\small $a', b'$ and $c'$ are consecutive non-critical edges of $P$.}
  \label{fig:rec1}
  \end{figure}

Let $P$ be the  $(\ge 5)$--gon of  $\ll$. Since $\ll\in \Im$ and $|\ll|=n$, $P$ is an $n$--gon. 
By  Observation~\ref{criticidad} we know that  $P$ has at most two critical edges. Hence, $P$ has
three consecutive non-critical edges, say $a', b'$ and $c'$. Let  $a,b$ and $c$ be, respectively, the pseudolines containing $a', b'$ and $c'$. Without loss of generality we  assume that  $a, b$ and $c$ are directed in such a way that $P$ lies to the left of each of them (see Figure~\ref{fig:rec1}). 

By Proposition~\ref{almenos-t} each of $a', b'$ and $c'$ is an edge of a triangle of  $\ll$. Let $A, B$ and $C$ be, respectively, the triangles of $\ll$ that are adjacent to $a', b'$ and $c'$. See Figure~\ref{fig:rec1}.

For $x\in \{a,b,c\}$ and $\ell=1,\ldots ,n-1$ let $x_{\ell}$ be the $\ell$--th pseudoline of  $\ll \setminus \{x\}$ that is crossed by $x$. Since $\ll$ is simple, the labels: $a_1,\ldots ,a_{n-1}$; $b_1,\ldots ,b_{n-1}$; and $c_1,\ldots ,c_{n-1}$ are well-defined. In fact, $\{a_1,\ldots , a_{n-1}\}\setminus \{c\}=\{c_1,\ldots , c_{n-1}\}\setminus \{a\}=\ll\setminus \{a,c\}$. From now on, we assume that $a_k=c$, $b_t=a$ and $c_r=a$. Then $i$) $A$ is formed by  $a, b=a_{k-1}$ and $a_{k-2}$, $ii$) $B$ is formed by  $b, a=b_t$ and $c=b_{t+1}$ and  $iii$) $C$ is formed by $c, b=c_{r+1}$ and $c_{r+2}$. See Figure~\ref{fig:rec1}.
 An easy consequence of  $i)$,  $ii$) and $iii)$ is, respectively, that  $3\le k\le n-1$,  $2\le t\le n-3$, and $1\le r\le n-3$.

\begin{figure}[ht]
  \begin{center}
  \includegraphics[width= 15cm, height=3cm]{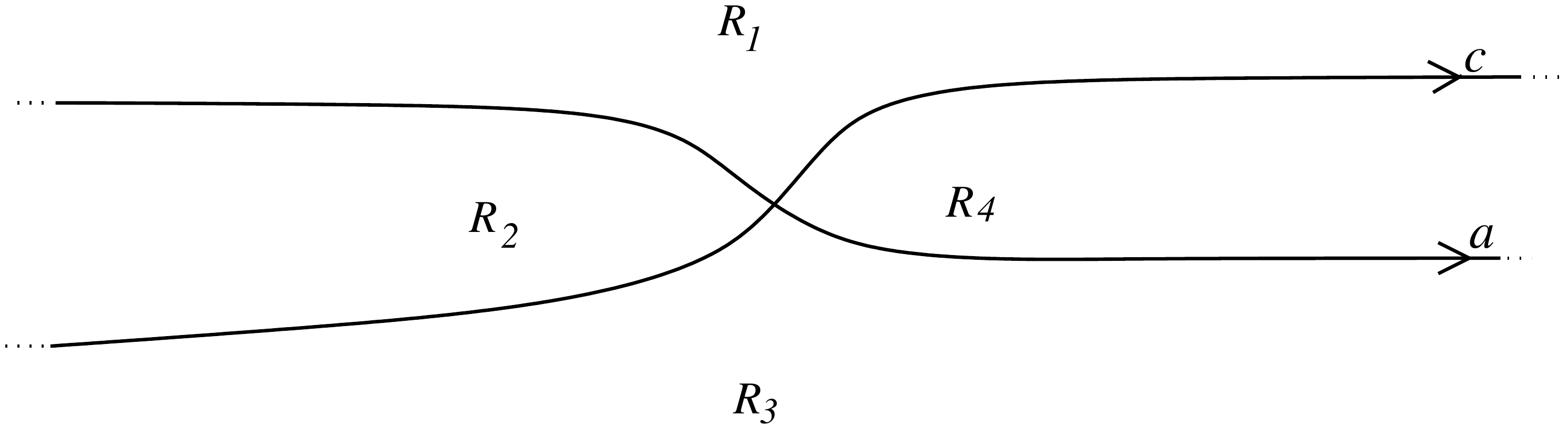}
  \end{center}
\caption{\small  $R_1, R_2, R_3$ and  $R_4$ are the four  semiplanes defined by $a$ and $c$.}
  \label{fig:rec2}
  \end{figure}

Let $R_1:=a^-\cap c^-$,  $R_2:=a^+\cap c^-$, $R_3:=a^+\cap c^+$ and  $R_4:=a^-\cap c^+$  be the four unbounded regions defined by  $a$ and $c$  (see Figure~\ref{fig:rec2}). Note that $P\subset R_1$.

\begin{figure}[ht]
  \begin{center}
  \includegraphics[width= 15cm, height=4cm]{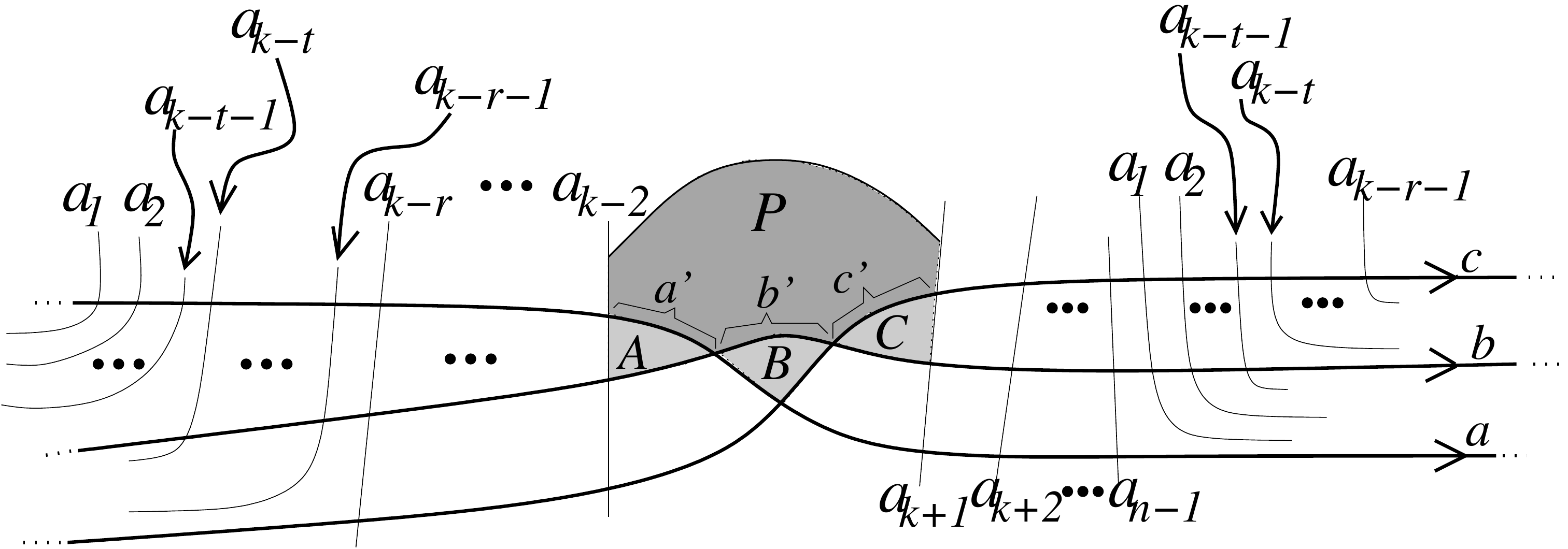}
  \end{center}
\caption{\small How $a_1,a_2,\ldots,a_{k-2},a_{k+1},a_{k+2},\ldots ,a_{n-1}$ are placed in  $R_2$ and $R_4$.}
  \label{fig:rec4}
  \end{figure}

Before we apply the inductive hypothesis, we need to establish some  of the structural properties of  $\ll$. Specifically, we shall prove that the pseudolines of $\ll\setminus \{a,b,c\}=\{a_1,\ldots,$ $a_{k-2}, a_{k+1},a_{k+2},\ldots ,a_{n-1}\}$ may be placed in $R_2$ and $R_4$ as shown in Figure~\ref{fig:rec4}. 

For $l',l''\in \{1,\ldots ,n-1\}$  and $w\in \{a,b,c\}$ let $w[l',l'']:=\{w_{l'}, w_{l'+1},\ldots ,w_{l''}\}$. By convention,  $w[l',l'']=\emptyset$ if $l'>l''$. Since $P$ has an edge in every pseudoline of $\ll$
 (since $\ll\in \Im$) and  
$P\subset R_1$, every pseudoline of  $\ll\setminus \{a,b,c\}$ crosses $R_1$. 

Since every  pseudoline of $\ll\setminus \{a,b,c\}$ intersects $R_1$ and crosses exactly once each of  $a, b$ and $c$,  $c[1,r-1]\subseteq b[1,t-1]\subseteq a[1,k-2]$  and $a[k+1,n-1]\subseteq b[t+2,n-1] \subseteq c[r+2,n-1]$. This implies that  $r\le t\le k-1$. 

Another consequence of the fact  that every  pseudoline of $\ll\setminus\{a,b,c\}$ is in $P\subset R_1$ is the following: if $x,y\in  a[1,k-2]$ or $x,y\in c[r+2,n-1]$ then $v_{x,y}\in R_1$. 

Because $R_2$ contains no crossings between pseudolines of $a[1,k-2]$,  and $c[1,r-1]\subseteq b[1,t-1]\subseteq a[1,k-2]$, 

\begin{equation}\label{eq:a1}
b_{t-j}=
\begin{cases}
a_{k-1-j}=c_{r-j} & \text{ if } j=1,2,\ldots ,r-1,\\
a_{k-1-j} & \text{ if } j=r,r+1,\ldots ,t-1. 
\end{cases}
\end{equation}

Analogously,  because $R_4$ contains no crossings  between pseudolines of $c[r+2,n-1]$, and  $a[k+1,n-1]\subseteq b[t+2,n-1] \subseteq c[r+2,n-1]$,

\begin{equation}\label{eq:c1}
b_{t+1+i}=
\begin{cases}
c_{r+1+i}=a_{k+i} & \text{ if } i=1,2,\ldots , n-k-1,\\
c_{r+1+i} & \text{ if } i=n-k,n-k+1,\ldots , n-t-2. 
\end{cases}
\end{equation}

\begin{figure}[ht]
  \begin{center}
  \includegraphics[width= 18cm, height=4.0cm]{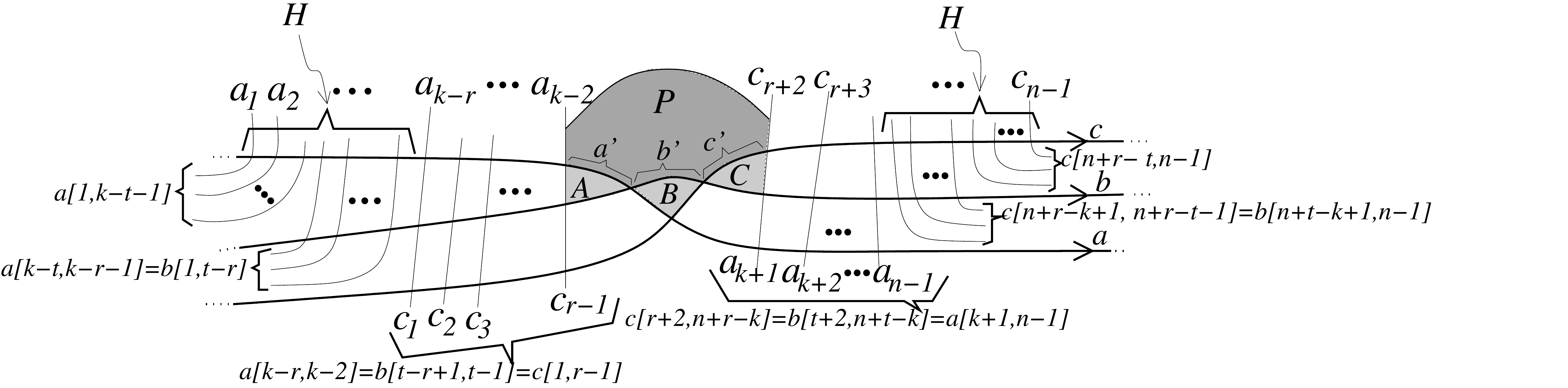}
  \end{center}
\caption{\small Since $R_2$ contains no crossings between pseudolines of $a[1,k-2]$, $\{a[1,k-t-1]$, $b[1,t-r]$,  $c[1,r-1]\}$ is a partition of $a[1,k-2]$. Analogously, 
$\{a[k+1,n-1], b[n+t-k+1,n-1], c[n+r-t,n-1]\}$ is a partition of $c[r+2,n-1]$   because $R_4$ contains no crossings between pseudolines of $c[r+2,n-1]$.}
  \label{fig:rec3}
  \end{figure}

It follows from  Equations  (\ref{eq:a1}) and (\ref{eq:c1}) that the pseudolines of 
$\ll\setminus\{a,b,c\}$ are placed in $R_2$ and $R_4$ as shown in Figure~\ref{fig:rec3}.

From Eqs.  (\ref{eq:a1}) and (\ref{eq:c1}) we deduce that $a[k-r,k-2]=c[1,r-1]$ and $a[k+1,n-1]=c[r+2,n+r-k]$, respectively. Let  $H:=a[1,k-r-1]$. From these equations and the fact that $a[1, k - 2]\cup  a[k + 1, n - 1] = c[1, r - 1] \cup c[r + 2, n - 1] =\ll \setminus\{a, b, c\} $  it follows that $H =c[n+r-k+1,n-1]$ (see Figure~\ref{fig:rec3}). Let $a_{i_1}=c_{j_1}$ and $a_{i_2}=c_{j_2}$ be pseudolines  of $H$. If $i_1<i_2$ then $j_1<j_2$  (respectively, if $i_1>i_2$ then $j_1>j_2$) because  $v_{a_{i_1},a_{i_2} }$ is in $R_1$. Thus, 
$c_{n+r-k+\ell}= a_{\ell}  \text{ for } \ell=1,\ldots ,k-r-1$. From this equation and Figure~\ref{fig:rec3} it follows that the pseudolines of $\ll\setminus \{a,b,c\} $  are placed in $R_2$ and $R_4$ as shown in Figure~\ref{fig:rec4}, as desired.



On the other hand, it is not difficult to see that the subarrangement of $\ll$ obtained by deleting $b$ belongs to  $\Im$.  By the induction hypothesis $\ll\setminus \{b\}$ is stretchable. Let $\ll^*_b$ be an arrangement of  straight lines, which is equivalent to $\ll\setminus \{b\}$.  If $\theta$ denotes an element of $\ll\setminus \{b\}$ (for example, a pseudoline, a region, etc.), we denote by  $\theta^*$ the corresponding element in $\ll^*_b$. For $s=1,\ldots , k-r-1$ let $m_s$ be the slope of $a^*_s$. Since any two lines of $\ll^*_b$ intersect exactly once,  $\ll^*_b$ has no lines with equal slopes, moreover, since the crossing of any  two lines of $H^*$ is in $R^*_1$, $m_{0}<m_1<\cdots <m_{k-r-1}<m_{k-r}$, where $m_{0}$ and $m_{k-r}$ are, respectively, the slopes of $a^*$ and $c^*$. Let $d^*$ be the line with slope  $(m_{k-t-1}+m_{k-t})/2$  through $v_{a^*,c^*}$. Since $\ll^*_b$ and $\ll\setminus \{b\}$ are equivalent and $m_{0}<m_1<\cdots <m_{k-r}$, $d^*$ crosses the lines of  $\ll^*_b\setminus \{a^*,c^*\}$ in the exact same order in which $b$ crosses the pseusolines of $\ll\setminus \{a,b,c\}$. Also, note that 
$d^*$ can be moved in such a way that (i) the order in which $d^*$ crosses the lines of  $\ll^*_b\setminus \{a^*,c^*\}$ is preserved  and (ii) $d^*$ intersects $R^*_1$ (conditions (i) and (ii) can be ensured by a sufficiently small motion, keeping the slope, of $d^*$). Finally, note that the arrangement of lines obtained by such a motion of $d^*$ is equivalent to $\ll$, as desired.  $\square$

 \begin{section}{On the number of non--isomorphic simple Euclidean arrangements with one $(\ge 5)$--gon.}\label{cuantos}

In this section we shall show that there are exponentially many non--isomorphic arrangements of 
$\Im$. We recall that  two Euclidean arrangements are isomorphic (or  {\it combinatorially equivalent}) if there is an incidence and dimension--preserving bijection between their induced cell decompositions
(see~\cite{felsner}).
 
 A $2$-colored necklace with $2m$ beads in which opposite beads have different colors is a {\it self-dual } necklace.  An example of a self-dual necklace is shown in Figure~\ref{fig:cont1} $i)$. Two self-dual necklaces are isomorphic if one can be obtained from the other by rotation or  reflection or both. In~\cite{palmer-robinson}  was proved that the number of non-isomorphic self-dual necklaces  is  given by 
 
 $$Q(m)=\frac{2^{\lfloor(m-1)/2 \rfloor}+\frac{1}{2m}\displaystyle\sum_{k|m,\hskip 0.1cm k\hskip 0.1cm odd} \phi(k)2^{m/k}}{2},$$
  
\noindent where $\phi(k)$ is the Euler totient function. 

Since $Q(m)$ grows exponentially with $m$, it is enough to exhibit a one-to-one correspondence 
between the set of self-dual necklaces and a subset of $\Im$. 

Let $C$ be a self-dual necklace with $2m\ge 6$  
beads colored 0 and 1, and let $P$ be the regular polygon of $2m$ sides. Now we  extend every edge of $P$ to both sides in such a way that each pair of non-parallel segments intersect as shown in  Figure~\ref{fig:cont1} $ii)$. Finally, we intersect every pair of parallel segments according to color  $1$ of $C$ as shown in Figure~\ref{fig:cont1} $iii)$. By construction, the resultant arrangement  $P_C$ belongs to  $\Im$. It is not difficult to see that (using this procedure) distinct necklaces generate distinct arrangements.

\begin{figure}[ht]
  \begin{center}
  \includegraphics[width= 14cm, height=5cm]{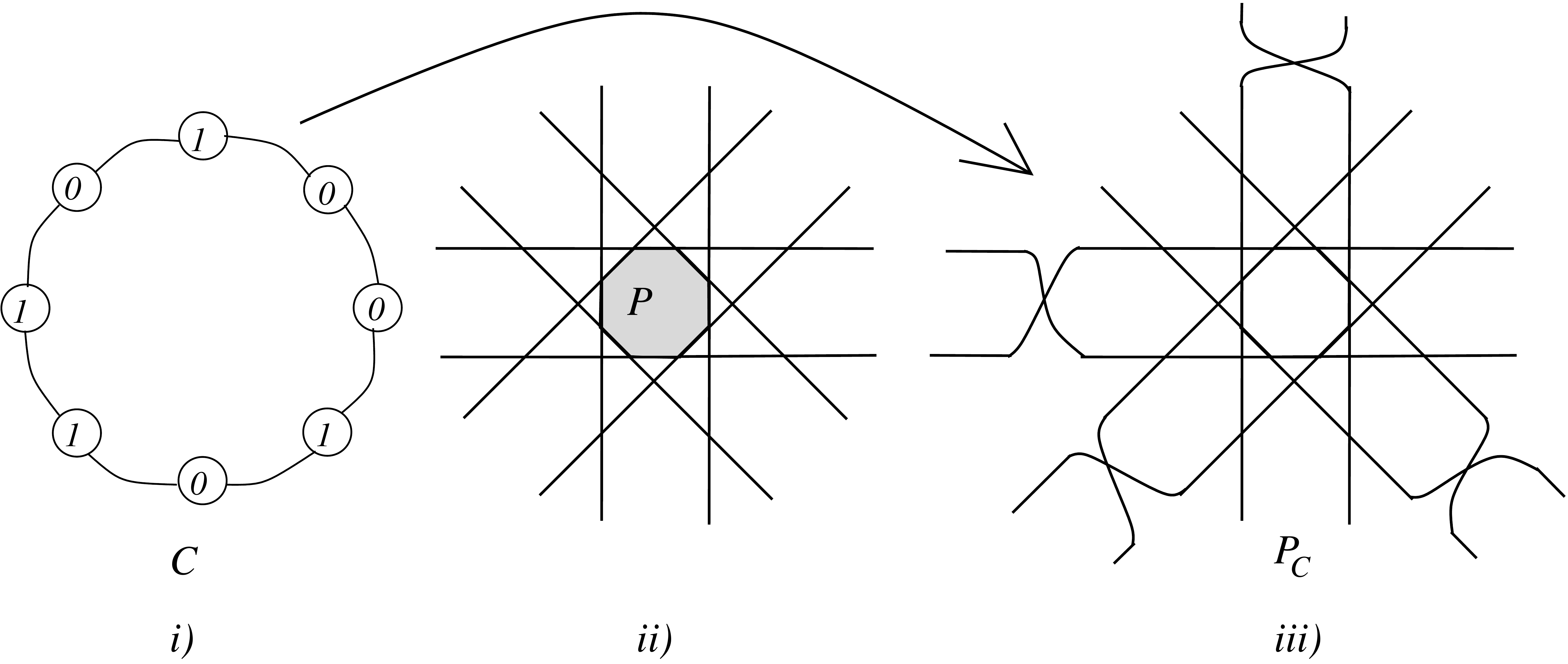}
  \end{center}
\caption{\small How to associate an arrangement of  $\Im$ to each self-dual necklace with $2m$ beads. $i)$ $C$ is a self-dual necklace with $8$ beads. $ii)$ The edges of the regular $8$--gon $P$ are  extended in such a way that any pair of non-parallel segments intersect. $iii)$ $P_C$ is the arrangement of  $\Im$ associated with $C$. We use the beads with color $1$ of $C$ to determine the side in that every pair of parallel segments intersect.}
  \label{fig:cont1}
  \end{figure}
 \end{section}

\end{document}